\documentclass{svmult}
\usepackage{bbm}
\usepackage{amsfonts}
\usepackage{mathrsfs}
\usepackage{amssymb}
\newenvironment{enumerate-roman}{\begin{enumerate}}{\end{enumerate}}
\usepackage{makeidx}         
\usepackage{graphicx}        
\usepackage{multicol}        
\usepackage[bottom]{footmisc}

\makeindex             

\newtheorem{thm}{Theorem}[section]
\newtheorem{lem}[thm]{Lemma}
\newtheorem{prop}[thm]{Proposition}

\newtheorem{rmk}[thm]{Remark}

\newtheorem{defi}[thm]{Definition}

\begin{document}

\baselineskip14pt

\title*{Stationary Stochastic Viscosity Solutions of SPDEs}
\titlerunning{Stationary Stochastic Viscosity Solutions of SPDEs}
\author{Qi Zhang\inst{\ {\rm 1}, {\rm 2}, {\rm 3}}}
\authorrunning{Q. Zhang}
\institute{$^{\rm 1}$ School of Mathematical Sciences, Fudan
University, Shanghai, 200433, China. (Current address)\\$^{\rm 2}$
School of
Mathematics, Loughborough University, Loughborough, LE11 3TU, UK.\\
$^{\rm 3}$ School of Mathematics, Shandong University, Jinan,
250100, China.\\\texttt{Email: qzh@fudan.edu.cn}}

\maketitle
\newcounter{bean}
\begin{abstract} In this paper we aim to find the stationary stochastic viscosity
solutions of a parabolic type SPDEs through the infinite horizon
backward doubly stochastic differential equations (BDSDEs). For
this, we study the existence, uniqueness and regularity of solutions
of the corresponding infinite horizon BDSDEs as well as the
``perfection procedure" applied to the solutions of BDSDEs. At last
the ``perfect"
stationary stochastic viscosity solutions of SPDEs constructed by solutions of corresponding BDSDEs are obtained.
\end{abstract}

\vskip5pt

\vskip15pt

\renewcommand{\theequation}{\arabic{section}.\arabic{equation}}

\section{Introduction}
\setcounter{equation}{0}

\ \ \ \ \ The pathwise stationary solution of a stochastic dynamical
system is one of the fundamental concept in the study of the long
time behaviour of the stochastic dynamical systems. It describes the
pathwise invariance of the stationary solution, over time, along the
measurable and $P$-preserving transformation
$\theta_t:\Omega\to\Omega$ and the pathwise limit of the solutions
of the random dynamical systems:
\begin{eqnarray}\label{zhao001}
u(t,Y(\omega),\omega)=Y(\theta_t\omega)\ \ \ t\geq 0,\ \rm{a.s.},
\end{eqnarray}
where $u: [0,\infty)\times U\times \Omega \to U$ is a measurable
random dynamical system on a measurable space $(U,\mathcal{B})$ over
a metric dynamical system ($\Omega$, $\cal F$, $P$,
$(\theta_t)_{t\geq0})$ and $Y:\Omega \to U$ is a $\cal F$-measurable
stationary solution. Needless to say that the ``one-force,
one-solution" setting is a natural extension of the equilibrium or
fixed point in the theory of the deterministic dynamical systems to
stochastic counterparts. Such a random fixed point consists of
infinitely many randomly moving invariant surfaces on the
configuration space due to the random external force pumped to the
system constantly. Therefore, in contrast to the deterministic
dynamical systems, the existence and stability of stationary
solutions of stochastic dynamical systems, generated e.g. by SDEs or
SPDEs, are a difficult and subtle problem.

In many works on random dynamical systems the existence of
stationary solutions is a basic assumption, e.g. in the study of
stability (Has$'$minskii \cite{ha}) and in the theory of stable and
unstable manifolds (Arnold \cite{ar}, Mohammed, Zhang and Zhao
\cite{mo-zh-zh}, Duan, Lu and Schmalfuss \cite{du-lu-sc1}). These
theories gave neither the existence of stationary solutions, nor a
way to find them. Although in \cite{mo-zh-zh}, Mohammed, Zhang and
Zhao introduced an integral equation of infinite horizon for the
stationary solutions of certain stochastic evolution equations, the
existence of the solutions of such stochastic integral equations in
general is far from clear.

Besides, from a pathwise stationary solution we can construct an
invariant measure for the skew product of the metric dynamical
system and the random dynamical system. The invariant measure
describes the invariance of a certain solution in law when time
changes, therefore it is a stationary measure of the Markov
transition probability. It is well known that an invariant measure
gives a stationary solution when it is a random Dirac measure.
Although an invariant measure of a random dynamical system on
${\mathbb{R}^{1}}$ gives a stationary solution, in general, this is
not true unless one considers an extended probability space.
However, considering the extended probability space, one essentially
regards the random dynamical system as noise as well, so the
dynamics is different. In fact, the pathwise stationary solution
gives the support of the corresponding invariant measure, so reveals
more detailed information than an invariant measure.

In spite of the importance of stationary solution, the difficulties,
arising mainly from random external force, prevent researchers form
finding a method universal to the stationary solutions of SPDEs with
great generalities. Some works on stationary solutions of certain
types of SPDEs usually under additive or linear noise include Sinai
\cite{si1}, \cite{si2} for stochastic Burgers' equations with
periodic or random forcing, Caraballo, Kloeden, Schmalfuss
\cite{kloeden} for stochastic evolution equations with small
Lipschitz constant. If one notices the solutions of infinite horizon
backward stochastic differential equations (BSDEs) give a classical
or viscosity solution of elliptic type PDEs (Poisson equations) from
the works of Peng \cite{pe} and Pardoux \cite{pa}, then it would be
natural to conjecture the stationary solutions of SPDEs can be
represented as the solutions of infinite horizon backward doubly
stochastic differential equations (BDSDEs). Inspired by this idea,
Zhang and Zhao in \cite{zh-zh1} proved that under the Lipschitz and
monotone conditions, the
$L_{\rho}^2({\mathbb{R}^{d}};{\mathbb{R}^{1}})\otimes
L_{\rho}^2({\mathbb{R}^{d}};{\mathbb{R}^{d}})$ valued solution of an
infinite horizon BDSDE exists and gives the stationary weak solution
of the corresponding parabolic SPDE. Zhang and Zhao
further considered this problem 
under the linear growth and monotone conditions in \cite{zh-zh2}. It
is easy to see the solutions of elliptic type PDEs give the
stationary solutions of the corresponding parabolic type PDEs,
however, for SPDEs of the parabolic type, such kind of connection
does not exist, so in this sense BDSDEs (or BSDEs) can be regarded
as more general SPDEs (or PDEs).

The stochastic viscosity solution of SPDE was first put forward by
Lions and Souganidis in \cite{li-so1} through stochastic
characteristics to remove the stochastic integrals in the SPDE. Then
Buckdahn and Ma in \cite{bu-ma1}-\cite{bu-ma3} gave their definition
through the Doss-Sussmann transformation. After that a few works on
stochastic viscosity solutions of SPDEs emerge using Buckdahn and
Ma's definition and corresponding BDSDEs, such as Boufoussi, Van
Casteren and Mrhardy \cite{bo-va-mr} for the SPDEs with Neumann
boundary conditions, Boufoussi and Mrhardy \cite{bo-mr} for the
multivalued SPDEs. Then an interesting question arises: can we also
find the stationary solution of some SPDE in the sense of stochastic
viscosity solution? This paper gives this question a positive
answer. By adopting Buckdahn and Ma's definition and using its
connection with BDSDE we can find the stationary stochastic
viscosity solution of the following SPDE:
\begin{eqnarray}\label{zz20}
v(t,x)&=&v(0,x)+\int_{0}^{t}[\mathscr{L}v(s,x)+f\big(x,v(s,x),\sigma^*(x)Dv(s,x)\big)]ds+\int_{0}^{t}\langle
g\big(x,v(s,x)\big),d{B}_s\rangle.
\end{eqnarray}
Here $(B_t)_{t\geq 0}$ is a Brownian motion with values in
$\mathbb{R}^l$; $f$, $g$ satisfy the condition (A.1)-(A.3) in
Section 2; $\mathscr{L}$ is the infinitesimal generator of the
diffusion process $X_{s}^{t,x}$ generated by the SDE as follows:
\begin{eqnarray}\label{qi17}
\left\{\begin{array}{l}
dX_{s}^{t,x}=b(X_{s}^{t,x})ds+\sigma(X_{s}^{t,x})dW_s,\ \ \ s>t\\
X_{s}^{t,x}=x,\ \ \ 0\leq s\leq t,
\end{array}\right.
\end{eqnarray}
where $(W_t)_{t\geq 0}$, independent of $(B_t)_{t\geq 0}$, is a
Brownian motion with values in $\mathbb{R}^d$ and $b$, $\sigma$
satisfy the condition (A.4) in Section 2, hence $\mathscr{L}$ is
given by
\begin{eqnarray*}\label{buchong1}
\mathscr{L}={1\over2}\sum_{i,j=1}^na_{ij}(x){{\partial^2}\over{\partial
x_i\partial x_j}}+\sum_{i=1}^nb_i(x){\partial\over{\partial x_i}}
\end{eqnarray*}
with $\big(a_{ij}(x)\big)=\sigma\sigma^*(x)$ and $(b_1(x),b_2(x),\cdot\cdot\cdot,b_n(x))^*=b(x)$. 
The infinite horizon BDSDEs we study as our tool can be written in
the following integration from:
\begin{eqnarray}\label{zz36}
{\rm e}^{-{K'\over2}s}Y_s^{t,x}&=&\int_{s}^{\infty}{\rm e}^{-{K'\over2}r}f(X_r^{t,x},Y_r^{t,x},Z_r^{t,x})dr+\int_{s}^{\infty}{K'\over2}{\rm e}^{-{K'\over2}r}Y_r^{t,x}dr\nonumber\\
&&-\int_{s}^{\infty}{\rm e}^{-{K'\over2}r}\langle
g(X_r^{t,x},Y_r^{t,x}),d^\dagger\hat{B}_r\rangle-\int_{s}^{\infty}{\rm
e}^{-{K'\over2}r}\langle Z_r^{t,x},dW_r\rangle,\ \ \ s\geq t.
\end{eqnarray}
Here 
$\hat{B}$ is the time reverse version of ${B}$, i.e. $\hat{B}_s=
B_{T-s}-B_T$ for arbitrary $T>0$ and all $s\in\mathbb{R}^{1}$,
and the integral w.r.t. $\hat{B}$ is a backward It$\hat {\rm o}$'s
integral (see \cite{zh-zh1} for details and the relationship between
the forward and backward It$\hat {\rm o}$'s integral). Our purpose
is to prove that, for arbitrary $T>0$ and $0\leq t\leq T$,
$v(t,x)(\omega)=Y_{T-t}^{T-t,x}(\hat{\omega})$ is a stationary
stochastic
viscosity solution of SPDE (\ref{zz20}). 
Five sections are organized in this paper for this purpose. In next
section we give brief introduction to the notion of stochastic
viscosity solutions of SPDEs and the connection between SPDEs and
BDSDEs in the sense of stochastic viscosity solution. In Section 3
under the assumption of the existence, uniqueness and regularity of
solution to infinite horizon BDSDE, we study its stationary
property, in which the general version ``perfection procedure" plays
an important role. The existence, uniqueness and regularity of
solution to infinite horizon BDSDE are proved in Section 4. In
Section 5 we deduce the stationary property for the stochastic
viscosity solutions of SPDEs constructed by the solutions of
infinite horizon BDSDEs.

As far as we know, the connection between the pathwise stationary
stochastic viscosity solutions of SPDEs and infinite horizon BDSDEs
in this paper is new. By the techniques as we dealt with the weak
solutions of PDEs or SPDEs in \cite{zh-zh2} and \cite{zh-zh3}, we
believe this connection can be extended to studying the stationary
stochastic viscosity solutions of more general parabolic SPDEs such
as those with linear or polynomial growth nonlinear terms, more
types of noises etc., but in this paper we only study Lipschitz
continuous nonlinear term and finite dimensional noise for
simplicity in order to initiate this method to the case of
stationary stochastic viscosity solutions of SPDEs. Finally we would
like to point out that the uniqueness of the stationary solution of
SPDE (\ref{zz20}) is still an open problem due to its complexity.

\section{Definition and Results for Stochastic Viscosity Solutions of SPDEs}\label{s25}
\setcounter{equation}{0}

\ \ \ \ The main purpose of this paper is to find the stationary
stochastic viscosity solution of SPDE (\ref{zz20}). As shown in
\cite{zh-zh1} and \cite{zh-zh2}, under appropriate conditions, for
$T\geq t\geq0$, defining $u(t,x)\triangleq v(T-t,x)$, we can obtain
the time reverse version of SPDE (\ref{zz20}) on $[0,T]$:
\begin{eqnarray}\label{zz21}
u(t,x)&=&u(T,x)+\int_{t}^{T}[\mathscr{L}u(s,x)+f\big(x,u(s,x),(\sigma^*\nabla
u)(s,x)\big)]ds-\int_{t}^{T}\langle
g\big(x,u(s,x)\big),d^\dagger\hat{B}_s\rangle.\nonumber\\
\end{eqnarray}
%
%
The BDSDE on $[t,T]$ associated with SPDE (\ref{zz21}) has the
following form:
\begin{eqnarray}\label{zz22}
Y_{s}^{t,x}=Y_{T}^{t,x}+\int_{s}^{T}f(X_{r}^{t,x},Y_{r}^{t,x},Z_{r}^{t,x})dr-\int_{s}^{T}\langle
g(X_{r}^{t,x},Y_{r}^{t,x}),d^\dagger\hat{B}_r\rangle-\int_{s}^{T}\langle
Z_{r}^{t,x},dW_r\rangle.
\end{eqnarray}

For $k,l\geq0$, we denote by $C_{b}^{k,l}$ 
the set of $C^{k,l}$-functions whose partial derivatives of order
for the first variable less than or equal to $k$ and for the second
variable less than or equal to $l$ are bounded. We assume
\begin{description}
\item[(A.1).] Functions $f: \mathbb{R}^d\times \mathbb{R}^1\times \mathbb{R}^{d}{\longrightarrow{\mathbb{R}^1}}$ and $g: \mathbb{R}^d\times \mathbb{R}^1{\longrightarrow{\mathbb{R}^{l}}}$ are $\mathscr{B}_{\mathbb{R}^{d}}\otimes\mathscr{B}_{\mathbb{R}^{1}}\otimes\mathscr{B}_{\mathbb{R}^{d}}$
and
$\mathscr{B}_{\mathbb{R}^{d}}\otimes\mathscr{B}_{\mathbb{R}^{1}}$
measurable respectively, and there exist constants $C_0$, $C_1$,
$C\geq0$ s.t. for any $(x_1, y_1, z_1)$, $(x_2, y_2, z_2)\in
\mathbb{R}^d\times \mathbb{R}^1\times \mathbb{R}^{d}$,
\begin{eqnarray*}
&&|f(x_1, y_1, z_1)-f(x_2, y_2, z_2)|^2\leq C_0|x_1-x_2|^2+C_1|y_1-y_2|^2+C|z_1-z_2|^2,\nonumber\\
&&|g(x_1, y_1)-g(x_2, y_2)|^2\leq C_0|x_1-x_2|^2+C|y_1-y_2|^2;
\end{eqnarray*}
\item[(A.2).] $g(\cdot,\cdot)\in
C_{b}^{2,3}(\mathbb{R}^d\times\mathbb{R}^1;\mathbb{R}^l)$;
\item[(A.3).] There exist constants $K\in\mathbb{R}^+$, $p>d+2$, $K<K'<2K$ and $\mu>0$ with $2\mu-{p\over2}K'-{p(p+1)\over2}C>0$ s.t. for
any $y_1$, $y_2\in \mathbb{R}^1$, $x$, $z\in \mathbb{R}^{d}$,
$$(y_1-y_2)(f(x, y_1, z)-f(x, y_2, z))\leq -\mu |y_1-y_2|^2;$$
\item[(A.4).] Functions $b(\cdot):
\mathbb{R}^d\longrightarrow \mathbb{R}^d$, $\sigma(\cdot):
\mathbb{R}^d\longrightarrow \mathbb{R}^{d\times d}$ are globally
Lipschitz continuous with Lipschitz constant $L$ and for $p$, $K$ in
(A.3), $K-pL-{p(p-1)\over2}L^2>0$.
\end{description}

Denote the set of $C^0$-functions with linear growth by $C^0_{l}$.
Buckdahn and Ma proved that if $u(T,\cdot)\in
C^0_{l}(\mathbb{R}^d;\mathbb{R}^1)$ is given, the solution
$Y_t^{t,x}$ of BDSDE (\ref{zz22}),
$(t,x)\in[0,T]\times\mathbb{R}^{d}$, is a stochastic viscosity
solution of SPDE (\ref{zz21}) under Conditions (A.1), (A.2) and
(A.4), therefore it gives the stochastic viscosity solution of SPDE
(\ref{zz20}) through the time reversal argument. To benefit the
reader, we include briefly Buckdahn and Ma's definition of
stochastic viscosity solution of SPDE (\ref{zz20}) through the
Doss-Sussmann transformation in \cite{bu-ma1}-\cite{bu-ma3}.

Let ${\cal N}$ be the class of $P$ null measure sets of
${\mathscr{F}}$. For any process $(\eta_t)_{t\geq0}$,
$\mathscr{F}_{s,t}^\eta\triangleq\sigma\{\eta_r-\eta_s$; ${0\leq
s\leq r\leq t}\}\bigvee{\cal N}$,
${\mathscr{F}}_t^\eta\triangleq{\mathscr{F}}_{0,t}^\eta$,
$\mathscr{F}_{t,\infty}^{\eta}\triangleq\bigvee_{T\geq0}{\mathscr{F}_{t,T}^\eta}$. Let $\mathbb{E}$ and $\mathbb{F}$ be the generic Euclidean spaces, then we denote\\
$\bullet$ $\mathscr{M}_{0,T}^B$ to be all the $\{\mathscr{F}_{t}^B\}_{t\geq0}$ stopping times $\tau$ such that $0\leq\tau\leq T$ a.s., where $T>0$ is some fixed time horizon;\\
$\bullet$ for any sub-$\sigma$-field $\mathscr{G}\subseteq\mathscr{F}_{T}^B$ and real number $p\geq0$, $L^p(\mathscr{G};\mathbb{E})$ to be $\mathbb{E}$-valued, $\mathscr{G}$-measurable random variables $\xi$ such that $E[|\xi|^p]<\infty$;\\
$\bullet$ for any sub-$\sigma$-field $\mathscr{G}\subseteq\mathscr{F}_{T}^B$, $C^{k,l}(\mathscr{G},[0,T]\times\mathbb{E};\mathbb{F})$ to be the space of all $C^{k,l}([0,T]\times\mathbb{E};\mathbb{F})$-valued random variables that are $\mathscr{G}\otimes\mathscr{B}_{[0,T]}\otimes\mathscr{B}_{\mathbb{E}}$-measurable;\\
$\bullet$
$C^{k,l}(\{\mathscr{F}_{t}^B\}_{t\geq0},[0,T]\times\mathbb{E};\mathbb{F})$
to be the space of all random fields $\varphi\in
C^{k,l}(\mathscr{F}_{T}^B,[0,T]\times\mathbb{E};\mathbb{F})$, such
that for fixed $x\in\mathbb{E}$, the mapping
$(t,\omega)\to\varphi(t,x,\omega)$ is
$\mathscr{F}_{t}^B$-progressively measurable.

The definition of stochastic viscosity solution depends heavily on
the following stochastic flow $\lambda\in
C^{0,0}(\{\mathscr{F}_{t}^B\}_{t\geq0},[0,T]\times\mathbb{R}^d\times\mathbb{R}^1;\mathbb{R}^1)$,
defined as the unique solution of the following SDE
\begin{eqnarray*}
\lambda(t,x,y)=y+{1\over2}\int_0^t\langle
g,D_yg\rangle(x,\lambda(s,x,y))ds-\int_0^t\langle
g(x,\lambda(s,x,y)),d{B}_s\rangle.
\end{eqnarray*}
Under Condition (A.2), $\lambda(t,x,y)$ is a stochastic flow, i.e.
for fixed $x$, the random field $\lambda(t,x,y)$ is continuously
differentiable in the variable $y$, and the mapping
$y\longrightarrow\lambda(t,x,y)$ defines a diffeomorphism for all
$(t,x)$, $P$-a.s. Denote the inverse of $\lambda$ by
$\zeta(t,x,y)=(\lambda(t,x,\cdot))^{-1}(y)$.
\begin{defi}\label{zz65} \rm{({\cite{bu-ma1}})}
A random field $w\in
C^{0,0}(\{\mathscr{F}_{t}^B\}_{t\geq0},[0,T]\times\mathbb{R}^d;\mathbb{R}^1)$
is called a stochastic viscosity subsolution (resp. supersolution)
of SPDE (\ref{zz20}), if $w(0,x)\leq$ (resp.$\geq$) $v(0,x)$,
$\forall x\in\mathbb{R}^d$; and if for any
$\tau\in\mathscr{M}_{0,T}^B$, $\xi\in
L^0(\mathscr{F}_{\tau}^B;\mathbb{R}^d)$, and any random field
$\varphi\in
C^{1,2}(\mathscr{F}_{\tau}^B,[0,T]\times\mathbb{R}^d;\mathbb{R}^1)$
satisfying
$$w(t,x)-\lambda(t,x,\varphi(t,x))\leq\ {\rm(resp.\ \geq)}\ 0=w(\tau,\xi)-\lambda(\tau,\xi,\varphi(\tau,\xi)),$$
for all $(t,x)$ in a neighborhood of $(\tau,\xi)$, P-a.e. on the set
$\{0<\tau<T\}$, it holds that
$$\mathscr{L}\psi(\tau,\xi)+f\big(\xi,\psi(\tau,\xi),\sigma^*(\xi)D\psi(\tau,\xi)\big)\geq\ {\rm(resp.\ \leq)}\ D_y\lambda\big(\tau,\xi,\varphi(\tau,\xi)\big)D_t\varphi(\tau,\xi),$$
P-a.e. on $\{0<\tau<T\}$, where
$\psi(t,x)\triangleq\lambda(t,x,\varphi(t,x))$.

A random field $w\in
C^{0,0}(\{\mathscr{F}_{t}^B\}_{t\geq0},[0,T]\times\mathbb{R}^d;\mathbb{R}^1)$
is called a stochastic viscosity solution of SPDE (\ref{zz20}), if
it is both a stochastic viscosity subsolution and a supersolution.
\end{defi}

By Doss-Sussmann transformation, SPDE (\ref{zz20}) can be converted
to the following PDE
\begin{eqnarray}\label{zz43}
\tilde{v}(t,x)=\tilde{v}(0,x)+\int_{0}^{t}[\mathscr{L}\tilde{v}(s,x)+\tilde{f}\big(s,x,\tilde{v}(s,x),\sigma^*(x)D\tilde{v}(s,x)\big)]ds,
\end{eqnarray}
where
\begin{eqnarray*}
&&\tilde{f}(t,x,y,z)={1\over
D_y\lambda(t,x,y)}\big(f(x,\lambda(t,x,y),\sigma^*(x)D_x\lambda(t,x,y)+D_y\lambda(t,x,y)z)\\
&&\ \ \ \ \ \ \ \ \ \ \ \ \ \ \ \ \ \ \ \ \ \ \ \ \ \ \ \ \ \ \ \ \
\ \
+\mathscr{L}_x\lambda(t,x,y)+\langle\sigma^*(x)D_{xy}\lambda(t,x,y),z\rangle+{1\over2}D_{yy}\lambda(t,x,y)|z|^2\big),
\end{eqnarray*}
and the stochastic viscosity solutions of (\ref{zz20}) and
(\ref{zz43}) have a kind of relationship like
$\tilde{v}(t,x)=\zeta(t,x,v(t,x))$. The Doss-Sussman transformation
plays a big role in the notion of the stochastic viscosity solution
of SPDE (\ref{zz20}). For more details, see Buckdahn and Ma
\cite{bu-ma1}-\cite{bu-ma3}.

Define
\begin{eqnarray*}
\mathscr{F}_{t,T}\triangleq{\mathscr{F}_{t,T}^{\hat{B}}}\bigvee
\mathscr{F}_t^W,\ {\rm for}\ 0\leq t\leq T;\ \ \
\mathscr{F}_t\triangleq{\mathscr{F}_{t,\infty}^{\hat{B}}}\bigvee
\mathscr{F}_t^W,\ {\rm for}\ t\geq0.
\end{eqnarray*}
For $q\geq2$, we define some useful solution spaces.
\begin{defi}\label{qi00}
Let $\mathbb{S}$ be a Banach space with norm $\|\cdot\|_\mathbb{S}$
and Borel $\sigma$-field $\mathscr{S}$. For $K\in\mathbb{R}^+$, we
denote by $M^{q,-K}([0,\infty);\mathbb{S})$ the set of
$\mathscr{B}_{\mathbb{R}^+}\otimes\mathscr{F}/\mathscr{S}$
measurable random processes $\{\phi(s)\}_{s\geq0}$ with values in
$\mathbb{S}$ satisfying
\begin{enumerate-roman}
\item $\phi(s):\Omega\rightarrow\mathbb{S}$ is $\mathscr{F}_s$ measurable for $s\geq 0$;
\item $E[\int_{0}^{\infty}{\rm e}^{-Ks}\|\phi(s)\|_\mathbb{S}^qds]<\infty$.
\end{enumerate-roman}
Also we denote by $S^{q,-K}([0,\infty);\mathbb{S})$ the set of
$\mathscr{B}_{\mathbb{R}^+}\otimes\mathscr{F}/\mathscr{S}$
measurable random processes $\{\psi(s)\}_{s\geq0}$ with values in
$\mathbb{S}$ satisfying
\begin{enumerate-roman}
\item $\psi(s):\Omega\rightarrow\mathbb{S}$ is $\mathscr{F}_s$ measurable for $s\geq0$ and $\psi(\cdot,\omega)$ is continuous $P$-a.s.;
\item $E[\sup_{s\geq0}{\rm
e}^{-Ks}\|\psi(s)\|_\mathbb{S}^q]<\infty$.
\end{enumerate-roman}
\end{defi}

Similarly, for $0\leq t\leq T<\infty$, we define
$M^{q,0}([t,T];\mathbb{S})$ and $S^{q,0}([t,T];\mathbb{S})$ on a
finite time interval.
\begin{defi}\label{zhao005}
Let $\mathbb{S}$ be a Banach space with norm $\|\cdot\|_\mathbb{S}$
and Borel $\sigma$-field $\mathscr{S}$. We denote by
$M^{q,0}([t,T];\mathbb{S})$ the set of
$\mathscr{B}_{[t,T]}\otimes\mathscr{F}/\mathscr{S}$ measurable
random processes $\{\phi(s)\}_{t\leq s\leq T}$ with values in
$\mathbb{S}$ satisfying
\begin{enumerate-roman}
\item $\phi(s):\Omega\rightarrow\mathbb{S}$ is $\mathscr{F}_{s,T}\bigvee{\mathscr{F}_{T,\infty}^{\hat{B}}}$ measurable for $t\leq s\leq T$;
\item $E[\int_{t}^{T}\|\phi(s)\|_\mathbb{S}^qds]<\infty$.
\end{enumerate-roman}
Also we denote by $S^{q,0}([t,T];\mathbb{S})$ the set of
$\mathscr{B}_{[t,T]}\otimes\mathscr{F}/\mathscr{S}$ measurable
random processes $\{\psi(s)\}_{t\leq s\leq T}$ with values in
$\mathbb{S}$ satisfying
\begin{enumerate-roman}
\item $\psi(s):\Omega\rightarrow\mathbb{S}$ is
$\mathscr{F}_{s,T}\bigvee{\mathscr{F}_{T,\infty}^{\hat{B}}}$
measurable for $t\leq s\leq T$ and $\psi(\cdot,\omega)$ is
continuous $P$-a.s.;
\item $E[\sup_{t\leq s\leq T}\|\psi(s)\|_\mathbb{S}^2]<\infty$.
\end{enumerate-roman}
\end{defi}

The following Buckdahn and Ma's result established the connection
between the solution of BDSDE (\ref{zz22}) and the stochastic
viscosity solution of SPDE (\ref{zz20}) on finite time interval
$[0,T]$.
\begin{thm}\label{theorem2.2} {\rm({\cite{bu-ma1}})}
Assume Conditions {\rm(A.1)}, {\rm(A.2)}, {\rm(A.4)} are satisfied
and the function $v(0,\cdot)\in C^0_{l}(\mathbb{R}^d)$ is given.
Then $v(t,x)=u(T-t,x)=Y_{T-t}^{T-t,x}$, where $Y^{t,x}_\cdot\in
S^{2,0}([0,T];\mathbb{R}^{1})$ is the solution of BDSDE
(\ref{zz22}), is a stochastic viscosity solution of SPDE
(\ref{zz20}) on finite time interval $[0,T]$.
\end{thm}
\begin{rmk}\label{zz24}
From the argument of Buckdahn and Ma we can see if we replace the
condition $v(0,\cdot)\in C^0_{l}(\mathbb{R}^d)$ in Theorem
\ref{theorem2.2} by that $v(0,x)$ is continuous w.r.t. $x$ and
$E[|v(0,X_T^{t,x})|^2]<\infty$, then the conclusion of Theorem
\ref{theorem2.2} remains true since
$E[|v(0,X_T^{t,x})|^2]=E[|Y_{T}^{t,x}|^2]<\infty$ 
guarantees the corresponding BDSDE has a square-integrable terminal
value.
\end{rmk}

\section{Stationary Property of Solutions of BDSDEs}
\setcounter{equation}{0}

\ \ \ \ The purpose of this section is to study the stationary
property of the solution to infinite horizon BDSDE (\ref{zz36}). In
order to show the main idea, we first assume that there exists a
unique solution $(Y^{t,x}_{\cdot}, Z^{t,x}_{\cdot})\in
S^{p,-K}([0,\infty); \mathbb{R}^1)\cap
M^{2,-K}([0,\infty);\mathbb{R}^1)\times M^{2,-K}([0,\infty);
\mathbb{R}^{d})$ to BDSDE (\ref{zz36}) and $(t,x)\to Y^{t,x}_{t}$ is
a.s. continuous. The study of the existence, uniqueness and
regularity of solution to BDSDE (\ref{zz36}) will be deferred to
next section.

We now construct the measurable metric dynamical system through
defining a measurable and measure-preserving shift. Let
$\hat{\theta}_t:\Omega\longrightarrow\Omega$, $t\geq0$, be a
measurable mapping on $(\Omega, {\mathscr{F}}, P)$, defined by
\begin{eqnarray*}
\hat{\theta}_{t}\circ \hat{B}_s=\hat{B}_{s+t}-\hat{B}_t,\ \ \
\hat{\theta}_{t}\circ W_s=W_{s+t}-W_t.
\end{eqnarray*}
Then for any $s,t\geq0$,
\begin{description}
\item[$(\textrm{i})$]$P\cdot\hat{\theta}_{t}^{-1}=P$;
\item[$(\textrm{i}\textrm{i})$]$\hat{\theta}_{0}=I$, where $I$ is the identity transformation on $\Omega$;
\item[$(\textrm{i}\textrm{i}\textrm{i})$]$\hat{\theta}_{s}\circ\hat{\theta}_{t}=\hat{\theta}_{s+t}$.
\end{description}
Also for an arbitrary $\mathscr{F}$ measurable random variable
$\phi$, set
\begin{eqnarray*}
\hat{\theta}\circ\phi(\omega)=\phi\big(\hat{\theta}(\omega)\big).
\end{eqnarray*}
%
For any $r\geq0$, $s\geq t$, $x\in\mathbb{R}^d$, applying
$\hat{\theta}_r$ to SDE (\ref{qi17}), we have
\begin{eqnarray*}
\hat{\theta}_r\circ
X_{s}^{t,x}=x+\int_{t+r}^{s+r}b(\hat{\theta}_r\circ
X_{u-r}^{t,x})du+\int_{t+r}^{s+r}\sigma(\hat{\theta}_r\circ
X_{u-r}^{t,x})dW_u.
\end{eqnarray*}
So under Condition (A.4), by the uniqueness of the solution, 
we have for any $r$, $t\geq0$, $x\in\mathbb{R}^{d}$,
\begin{eqnarray}\label{qi18}
\hat{\theta}_r\circ X_{s}^{t,x}=X_{s+r}^{t+r,x},\ \ {\rm for}\ {\rm
all}\ s\geq0\ {\rm a.s.}
\end{eqnarray}
For $Y\in\mathbb{R}^{1}$, $x$, $Z\in\mathbb{R}^{d}$, let
\begin{eqnarray*}
\hat{f}(\mathcal{T},Y,Z)=f(X_{s}^{t,x},Y,Z),\ \ \
\hat{g}(\mathcal{T},Y,Z)=g(X_{s}^{t,x},Y,Z).
\end{eqnarray*}
Here we take $\mathcal{T}=(s,t)$ as a dual time variable (t is
fixed). Using (\ref{qi18}) we can verify that $\hat{f}$ and
$\hat{g}$ satisfy the stationary conditions in Proposition 2.5 in
\cite{zh-zh1} for any $\hat{\theta}_{r}$ $(r\geq0)$, $\mathcal{T}$,
$Y$ and $Z$, then using a similar argument as in Theorem 2.12 in
\cite{zh-zh1} we can deduce the following proposition by the
uniqueness of BDSDE (\ref{zz36}):
\begin{prop}\label{qi031}
Assume BDSDE (\ref{zz36}) has a unique solution $(Y^{t,x}_{\cdot},
Z^{t,x}_{\cdot})\in S^{p,-K}([0,\infty); \mathbb{R}^1)\cap
M^{2,-K}([0,\infty);\mathbb{R}^1)\times M^{2,-K}([0,\infty);
\mathbb{R}^{d})$, then under Condition (A.4), $(Y^{t,x}_s
Z^{t,x}_s)_{s\geq0}$ satisfies the following stationary property
w.r.t. $\hat{\theta}_\cdot$: for any $r$, $t\geq0$,
$x\in\mathbb{R}^{d}$,
\begin{eqnarray*}
\hat{\theta}_r\circ Y^{t,x}_s=Y^{t+r,x}_{s+r}, \ \
\hat{\theta}_r\circ Z^{t,x}_s=Z^{t+r,x}_{s+r}\ \ {\rm for}\ {\rm
all}\ s\geq0\ {\rm a.s.}
\end{eqnarray*}
In particular, for any $r$, $t\geq0$, $x\in\mathbb{R}^{d}$,
\begin{eqnarray}\label{zz18}
\hat{\theta}_r\circ Y^{t,x}_t=Y^{t+r,x}_{t+r}\ \ \ {\rm a.s.}
\end{eqnarray}
\end{prop}

If we regard $Y_t^{t,x}$ as a function of $(t,x)$, (\ref{zz18})
gives a ``very crude" stationary property of $Y$. Borrowing the idea
of perfecting crude cocycles in \cite{ar} and \cite{ar-sc}, we then
prove the following theorem which makes the ``very crude" stationary
property
of $Y$ ``perfect". 
\begin{thm}\label{qi032}
Let $(\Omega,\mathscr{F},P)$ be a probability space and $\mathbb{H}$
be a separable Hausdorff topological space with $\sigma$-algebra
$\mathscr{H}$. Assume $Y(t,x,\omega)$:
$[0,\infty)\times\mathbb{R}^d\times\Omega\longrightarrow \mathbb{H}$
is
$\mathcal{B}_{\mathbb{R}^+}\otimes\mathcal{B}_{\mathbb{R}^d}\otimes{\mathscr{F}}$
measurable, a.s. continuous w.r.t. $t$, $x$ and satisfies the ``very
crude" stationary property w.r.t. $\hat{\theta}_\cdot$, i.e. for any
$t,r\geq0$, $x\in\mathbb{R}^d$
\begin{eqnarray}\label{zhang100}
\hat{\theta}_{r}\circ Y(t,x,\omega)=Y(t+r,x,\omega)\ \ \rm{a.s.}
\end{eqnarray}
Then there exists a $\hat{Y}(t,x,\omega)$ which is an
indistinguishable version of ${Y}(t,x,\omega)$ s.t.
$\hat{Y}(t,x,\omega)$ is
$\mathcal{B}_{\mathbb{R}^+}\otimes\mathcal{B}_{\mathbb{R}^d}\otimes{\mathscr{F}}$
measurable, continuous w.r.t. $t$, $x$ for all $\omega$ and
satisfies the ``perfect" stationary property w.r.t.
$\hat{\theta}_\cdot$:
\begin{eqnarray}\label{zz63}
\hat{\theta}_{r}\circ \hat{Y}(t,x,\omega)=\hat{Y}(t+r,x,\omega)\ \ \
{\rm for}\ {\rm all}\ t,r\geq0,\ x\in\mathbb{R}^d\ {\rm a.s.}
\end{eqnarray}
\end{thm}
{\em Proof}. From the continuity of $Y(t,x,\omega)$ w.r.t. $t$, $x$
and using a standard argument, we easily see that for any $r\geq0$,
\begin{eqnarray}\label{zhang331}
\hat{\theta}_{r}\circ Y(t,x,\omega)=Y(t+r,x,\omega)\ \ \ {\rm for}\
{\rm all}\ t\geq0,\ x\in\mathbb{R}^d\ {\rm a.s.}
\end{eqnarray}
Define
\begin{eqnarray*}
&&M=\{(r,\omega):\hat{\theta}_{r}\circ Y(t,x,\omega)=Y(t+r,x,\omega)\ {\rm for}\ {\rm all}\ t,x\};\\
&&\tilde{\Omega}=\{\omega:(r,\omega)\in M\ {\rm for}\ {\rm a.e.}\ r\};\\
&&{\Omega^*}=\{\omega:\hat{\theta}_r\omega\in\tilde{\Omega}\ {\rm for}\ {\rm a.e.}\ r\};\\
&&A(r,t,x,\omega)=\hat{\theta}_r\circ Y(t,x,\omega)-Y(t+r,x,\omega).
\end{eqnarray*}
Obviously, $A(r,t,x,\omega)$ is measurable w.r.t.
$\mathcal{B}_{\mathbb{R}^+}\otimes\mathcal{B}_{\mathbb{R}^+}\otimes\mathcal{B}_{\mathbb{R}^d}\otimes{\mathscr{F}}$.
If we denote by $Q$ and $\tilde{Q}$ the normalized Lebesgue measure
on $\mathbb{R}^+$ and $\mathbb{R}^d$ respectively such that
$Q(\mathbb{R}^+)=1$ and $\tilde{Q}(\mathbb{R}^d)=1$, then by
(\ref{zhang331}),
\begin{eqnarray}\label{zhang351}
Q\otimes Q\otimes\tilde{Q}\otimes
P\big(A^{-1}(0)\big)=\int_{\mathbb{R}^+}\int_{\mathbb{R}^+}\int_{\mathbb{R}^d}\int_{\Omega}I_{A^{-1}(0)}(r,t,x,\omega)dPd\tilde{Q}dQdQ=1,
\end{eqnarray}
where $I$ is the indicator function in
$(\mathbb{R}^+\times\mathbb{R}^+\times\mathbb{R}^d\times\Omega,\
\mathcal{B}_{\mathbb{R}^+}\otimes\mathcal{B}_{\mathbb{R}^+}\otimes\mathcal{B}_{\mathbb{R}^d}\otimes\mathscr{F})$.
It is easy to see that
\begin{eqnarray*}
M=\{(r,\omega):\int_{\mathbb{R}^+}\int_{\mathbb{R}^d}I_{A^{-1}(0)}(r,t,x,\omega)d\tilde{Q}dQ=1\}\in\mathcal{B}_{\mathbb{R}^+}\otimes\mathscr{F}.
\end{eqnarray*}
And by (\ref{zhang351}), we have
\begin{eqnarray*}
Q\otimes P(M)=Q\otimes
P\big(\{(r,\omega):\int_{\mathbb{R}^+}\int_{\mathbb{R}^d}I_{A^{-1}(0)}(r,t,x,\omega)d\tilde{Q}dQ=1\}\big)=1.
\end{eqnarray*}
Similarly, we can also know
\begin{eqnarray*}
\tilde{\Omega}=\{\omega:\int_{\mathbb{R}^+}I_M(r,\omega)dQ=1\}\in\mathscr{F}
\end{eqnarray*}
and
\begin{eqnarray*}
P(\tilde{\Omega})=P\big(\{\omega:\int_{\mathbb{R}^+}I_M(r,\omega)dQ=1\}\big)=1.
\end{eqnarray*}
Moreover, the measurability of $\Omega^*$ can be seen easily as
\begin{eqnarray*}
{\Omega^*}=\{\omega:\int_{\mathbb{R}^+}\int_{\mathbb{R}^+}I_M(r,\hat{\theta}_u\omega)dQdQ=1\}\in\mathscr{F}.
\end{eqnarray*}
And since $\tilde{\Omega}$ has full measure,
\begin{eqnarray*}
P({\Omega^*})&\geq&P\big(\{\omega:Y(t+r,x,\hat{\theta}_u\omega)=Y(t,x,\hat{\theta}_r\circ\hat{\theta}_u\omega)\ {\rm for}\ {\rm a.e.}\ r\ {\rm and}\ u,\ {\rm and}\ {\rm all}\ t,\ x\}\bigcap{\tilde{\Omega}}\big)\\
&=&P\big(\{\omega:Y(t+r+u,x,\omega)=Y(t,x,\hat{\theta}_{r+u}\omega)\ {\rm for}\ {\rm a.e.}\ r\ {\rm and}\ u,\ {\rm and}\ {\rm all}\ t,\ x\}\bigcap{\tilde{\Omega}}\big)\\
&=&P\big(\{\omega:Y(t+r',x,\omega)=Y(t,x,\hat{\theta}_{r'}\omega)\ {\rm for}\ {\rm a.e.}\ r',\ {\rm and}\ {\rm all}\ t,\ x\}\bigcap{\tilde{\Omega}}\big)\\
&=&P(\tilde{\Omega})\\
&=&1.
\end{eqnarray*}
One can prove $\hat{\theta}_u\Omega^*\subset\Omega^*$ for any
$u\geq0$. Indeed, for any $\omega\in\hat{\theta}_u\Omega^*$, there
exists $\hat{\omega}\in\Omega^*$ s.t.
$\omega=\hat{\theta}_u\hat{\omega}$ and
$\hat{\theta}_r\hat{\omega}\in\tilde{\Omega}$ for a.e. $r\geq0$. But
$\hat{\theta}_r\omega=\hat{\theta}_{u+r}\hat{\omega}\in\tilde{\Omega}$
for a.e. $r\geq0$, so $\omega\in\Omega^*$. That is to say
$\hat{\theta}_u\Omega^*\subset\Omega^*$. Define
\begin{eqnarray*}
\left\{\begin{array}{l}
\hat{Y}(t,x,\omega)=Y(t-r,x,\hat{\theta}_r\omega),\ \ \ {\rm where}\ r\in[0,t]\ {\rm with}\ \hat{\theta}_r\omega\in\tilde{\Omega},\ {\rm if}\ \omega\in\Omega^*,\nonumber\\
\hat{Y}(t,x,\omega)=0,\ \ \ \ \ \ \ \ \ \ \ \ \ \ \ \ \ \ \ \ \ \ \
\ \ \ \ \ \ \ \ \ \ \ \ \ \ \ \ \ \ \ \ \ \ \ \ \ \ \ \ \ \ \ {\rm
if}\ \omega\in{\Omega^*}^c.\nonumber
\end{array}\right.
\end{eqnarray*}
An important fact is that if $\omega\in\Omega^*$, then for an
arbitrary $r\in[0,t]$ with $\hat{\theta}_r\omega\in\tilde{\Omega}$,
$Y(t-r,x,\hat{\theta}_r\omega)$ is independent of $r$ and
\begin{eqnarray}\label{qi11}
{Y}(t-r,x,\hat{\theta}_{r}\omega)=Y(t,x,\omega).
\end{eqnarray}
To see this, as $\hat{\theta}_r\omega\in\tilde{\Omega}$, so there
exists $u\geq r$ s.t. $(u,\hat{\theta}_r\omega)\in M$ and
$(u-r,\hat{\theta}_r\omega)\in M$. If not, it means for a.e. $r$
there doesn't exist $u$ satisfying $(u,\hat{\theta}_r\omega)\in M$
and $(u-r,\hat{\theta}_r\omega)\in M$. Then one can easily get the
measure of $\{u:(u,\hat{\theta}_r\omega)\notin M\}$ is positive.
That is a contradiction. So such a $u$ certainly exists and
satisfies
\begin{eqnarray*}
\hat{\theta}_u{Y}(t-r,x,\hat{\theta}_{r}\omega)=Y(t-r+u,x,\hat{\theta}_r\omega)=Y(t,x,\hat{\theta}_{u-r}\hat{\theta}_r\omega)=Y(t,x,\hat{\theta}_u\omega).
\end{eqnarray*}
So
\begin{eqnarray*}
{Y}(t-r,x,\hat{\theta}_{r}\omega)=\hat{\theta}_u^{-1}Y(t,x,\hat{\theta}_u\omega)=Y(t,x,\omega).
\end{eqnarray*}
Therefore (\ref{qi11}) is true and $\hat{Y}(t,x,\omega)$ doesn't
depend on the choice of $r$. That is to say $\hat{Y}(t,x,\omega)$ is
well defined. Moreover (\ref{qi11}) implies that
${Y}(t,x,\omega)=\hat{Y}(t,x,\omega)$ for all $t\geq0$,
$x\in\mathbb{R}^d$ on a full measure set $\Omega^*$, thus
${Y}(t,x,\omega)$ and $\hat{Y}(t,x,\omega)$ are indistinguishable.
Define
\begin{eqnarray*}
\left\{\begin{array}{l}
B(r,t,x,\omega)=Y(t-r,x,\hat{\theta}_r\omega),\ \ \ \ \ \ {\rm if}\ r\in[0,t],\ \ \hat{\theta}_r\omega\in\tilde{\Omega},\ \ {\rm and}\ \omega\in\Omega^*,\nonumber\\
B(r,t,x,\omega)=0,\ \ \ \ \ \ \ \ \ \ \ \ \ \ \ \ \ {\rm
otherwise}.\nonumber
\end{array}\right.
\end{eqnarray*}
Then $B(r,t,x,\omega)$ is
$\mathcal{B}_{\mathbb{R}^+}\otimes\mathcal{B}_{\mathbb{R}^+}\otimes\mathcal{B}_{\mathbb{R}^d}\otimes\mathscr{F}$
measurable. By the definition of $\Omega^*$, if $\omega\in\Omega^*$,
then for a.e. $0\leq r\leq t$,
$\hat{\theta}_r\omega\in\tilde{\Omega}$. We denote $L(r)$ the
Lebesgue measure in $[0,t]$. Since the countable base of $H$
generates $\mathscr{H}$ and separates points, $(H,\mathscr{H})$ is
isomorphic as a measurable space to a subset of $[0,1]$.
Consequently, for all $t,x,\omega$,
\begin{eqnarray*}
\hat{Y}(t,x,\omega)=\int_{0}^{t}B(r,t,x,\omega)dL(r).
\end{eqnarray*}
So by Fubini's theorem, $\hat{Y}(t,x,\omega)$ is
$\mathcal{B}_{\mathbb{R}^+}\otimes\mathcal{B}_{\mathbb{R}^d}\otimes\mathscr{F}$
measurable. $\hat{Y}(t,x,\omega)$ is a.s continuous w.r.t. $t$, $x$
due to the a.s continuity of ${Y}(t-r,x,\omega)$. But there exists a
null measure set $N\in\mathscr{F}$ s.t. $\{\omega:\
\hat{Y}(t,x,\omega)\ {\rm is}\ {\rm not}\ {\rm continuous}\ {\rm
w.r.t.}\ t,x\}\subset N$. Let $\hat{Y}(t,x,\omega)$ on $N$ equal
$0$. We still denote this new version of $\hat{Y}(t,x,\omega)$ by
$\hat{Y}(t,x,\omega)$, then this version of $\hat{Y}(t,x,\omega)$ is
continuous for all $\omega$.

The remaining work is to check $\hat{Y}(t,x,\omega)$ satisfies the
``perfect" stationary property (\ref{zz63}). For $\omega\in\Omega^*$
and all $r\geq0$,
$\hat{\theta}_r\omega\in\hat{\theta}_r\Omega^*\subset\Omega^*$. Pick
a $u$ s.t. $\hat{\theta}_u\omega\in\tilde{\Omega}$,
$\hat{\theta}_{u+r}\omega\in\tilde{\Omega}$, then by (\ref{qi11}) we
have
\begin{eqnarray*}
\hat{Y}(t,x,\hat{\theta}_{r}\omega)&=&Y(t-u,x,\hat{\theta}_{u+r}\omega)=Y(t+r-u-r,x,\hat{\theta}_{u+r}\omega)\\
&=&Y(t+r,x,\omega)={Y}(t+r-u,x,\hat{\theta}_u\omega)=\hat{Y}(t+r,x,\omega).
\end{eqnarray*}
The theorem is proved.$\hfill\diamond$


From now on, we neglect the difference between two distinguishable
random processes. Then with Proposition \ref{qi031} and Theorem
\ref{qi032}, it follows immediately that
\begin{thm}\label{zz42}
If BDSDE (\ref{zz36}) has a unique solution $(Y^{t,x}_{\cdot},
Z^{t,x}_{\cdot})\in S^{p,-K}([0, \infty);\mathbb{R}^1)\cap
M^{2,-K}([0,\infty);\mathbb{R}^1)\\\times M^{2,-K}([0, \infty);
\mathbb{R}^{d})$ and $(t,x)\to Y_t^{t,x}$ is a.s. continuous, then
under Condition (A.4), $Y_t^{t,x}$ satisfies the ``perfect"
stationary property w.r.t. $\hat{\theta}_\cdot$, i.e.
\begin{eqnarray*}
\hat{\theta}_r\circ Y_t^{t,x}=Y_{t+r}^{t+r,x}\ \ {\rm for}\ {\rm
all}\ r,\ t\geq0,\ x\in\mathbb{R}^d\ {\rm a.s.}
\end{eqnarray*}
\end{thm}


\section{Infinite Horizon BDSDEs}\label{s24}
\setcounter{equation}{0}

\ \ \ \ In this section we first prove the assumption in Theorem
\ref{zz42} that BDSDE (\ref{zz36}) has a unique solution
$(Y_{\cdot}, Z_{\cdot})\in S^{p,-K}\bigcap
M^{2,-K}([0,\infty);\mathbb{R}^1)\times M^{2,-K}([0, \infty);
\mathbb{R}^{d})$ is obtainable and reasonable under Conditions
(A.1)-(A.4).
To begin with, we briefly introduce the pioneering work by Pardoux
and Peng in \cite{pa-pe3} for the following finite horizon BDSDE:
\begin{eqnarray}\label{zz27}
Y_s=Y_T+\int_{s}^{T}f(r,Y_r,Z_r)dr-\int_{s}^{T}\langle
g(r,Y_r,Z_r),d^\dagger\hat{B}_r\rangle-\int_{s}^{T}\langle
Z_r,dW_r\rangle.
\end{eqnarray}
Here we only consider $\mathbb{R}^1$-valued BDSDE for our purpose.
One can also refer to \cite{pa-pe3} for multi-dimensional BDSDE if
interested. Assume
\begin{description}
\item[(A.1)$'$.] Functions $f: \Omega\times[0,T]\times \mathbb{R}^1\times \mathbb{R}^{d}{\longrightarrow{\mathbb{R}^1}}$ and $g: \Omega\times[0,T]\times \mathbb{R}^1\times \mathbb{R}^{d}{\longrightarrow{\mathbb{R}^{l}}}$ are ${\mathscr{F}}\otimes\mathscr{B}_{[0,T]}\otimes\mathscr{B}_{\mathbb{R}^{1}}\otimes\mathscr{B}_{\mathbb{R}^{d}}$
measurable, and for any $(y, z)\in\mathbb{R}^1\times\mathbb{R}^{d}$,
$f(\cdot,y,z)\in M^{2,0}([0,T];\mathbb{R}^{1})$ and $g(\cdot,y,z)\in
M^{2,0}([0,T];\mathbb{R}^{d})$, moreover there exist constants
$C\geq0$ and $0\leq\alpha<1$ s.t. for any $r\in[0,T]$, $(y_1, z_1)$,
$(y_2, z_2)\in\mathbb{R}^1\times \mathbb{R}^{d}$,
\begin{eqnarray*}
&&|f(r, y_1, z_1)-f(r, y_2, z_2)|^2\leq C|y_1-y_2|^2+C|z_1-z_2|^2,\nonumber\\
&&|g(r, y_1, z_1)-g(r, y_2, z_2)|^2\leq
C|y_1-y_2|^2+\alpha|z_1-z_2|^2.
\end{eqnarray*}
\end{description}
\begin{thm}\label{theorem2.1}
{\rm(\cite{pa-pe3})} Under Condition {\rm(A.1)$'$}, for any given
${\cal F}_T\bigvee{\mathscr{F}_{T,\infty}^{\hat{B}}}$ measurable
$Y_T\in L^2(\Omega)$, BDSDE (\ref{zz27}) has a unique solution
\begin{center}
$(Y_{\cdot}, Z_{\cdot})\in S^{2,0}([0,T];\mathbb{R}^1)\bigotimes
M^{2,0}([0,T]; \mathbb{R}^{d})$.
\end{center}
\end{thm}
In \cite{pa-pe3}, Pardoux and Peng also discussed a type of forward
BDSDE, a special case of BDSDE (\ref{zz27}),
\begin{eqnarray}\label{zz28}
Y_s^{t,x}&=&h(X_T^{t,x})+\int_s^Tf(X_r^{t,x},Y_r^{t,x},Z_r^{t,x})dr\nonumber\\
&&-\int_s^T\langle
g(X_r^{t,x},Y_r^{t,x},Z_r^{t,x}),d^\dagger\hat{B}_r\rangle-\int_s^T\langle
Z_r^{t,x},dW_r\rangle,
\end{eqnarray}
where $(X_s^{t,x})_{t\leq s\leq T}$ is the solution of SDE
(\ref{qi17}). Assume
\begin{description}
\item[(A.2)$'$.] Functions $f: \mathbb{R}^{d}\times\mathbb{R}^1\times\mathbb{R}^{d}{\longrightarrow{\mathbb{R}^1}}$ and $g: \mathbb{R}^{d}\times\mathbb{R}^1\times\mathbb{R}^{d}{\longrightarrow{\mathbb{R}^{l}}}$ are $\mathscr{B}_{\mathbb{R}^{d}}\otimes\mathscr{B}_{\mathbb{R}^{1}}\otimes\mathscr{B}_{\mathbb{R}^{d}}$
measurable, and there exist constants $C\geq0$ and $0\leq\alpha<1$
s.t. for any $(x_1, y_1, z_1)$, $(x_2, y_2, z_2)\in
\mathbb{R}^d\times \mathbb{R}^1\times \mathbb{R}^{d}$,
\begin{eqnarray*}
&&|f(x_1, y_1, z_1)-f(x_2, y_2, z_2)|^2\leq C|x_1-x_2|^2+C|y_1-y_2|^2+C|z_1-z_2|^2,\nonumber\\
&&|g(x_1, y_1, z_1)-g(x_2, y_2, z_2)|^2\leq
C|x_1-x_2|^2+C|y_1-y_2|^2+\alpha|z_1-z_2|^2.
\end{eqnarray*}
\end{description}
For BDSDE (\ref{zz28}), it is not difficult to deduce from Theorem
\ref{theorem2.1} that
\begin{thm}\label{zz38} Under Condition {\rm(A.2)$'$}, for each $x\in\mathbb{R}^{d}$ and any
given ${\cal F}_T\bigvee{\mathscr{F}_{T,\infty}^{\hat{B}}}$
measurable $h$ satisfying $h(X_T^{t,x})\in L^2(\Omega)$, BDSDE
(\ref{zz28}) has a unique solution
\begin{center}
$(Y^{t,x}_{\cdot}, Z^{t,x}_{\cdot})\in
S^{2,0}([t,T];\mathbb{R}^1)\bigotimes M^{2,0}([t,T];
\mathbb{R}^{d})$.
\end{center}
\end{thm}

In \cite{pa-pe3}, for the first time, Pardoux and Peng associated
the classical solution of SPDE, if any, with the solution of BDSDE
(\ref{zz28}). They proved that under some strong smoothness
conditions of $h$, $b$, $\sigma$, $f$ and $g$ (for details see
\cite{pa-pe3}), $u(t,x)=Y_t^{t,x}$, where $Y$ is the unique solution
of BDSDE (\ref{zz28}), $(t,x)\in[0,T]\times\mathbb{R}^{d}$, is
independent of ${{\cal F}_T^W}$ and is the unique classical solution
of the following backward SPDE
\begin{eqnarray*}\label{zz29}
u(t,x)&=&h(x)+\int_{t}^{T}[\mathscr{L}u(s,x)+f\big(x,u(s,x),\sigma^*(x)Du(s,x)\big)]ds\nonumber\\
&&-\int_{t}^{T}\langle
g\big(x,u(s,x),\sigma^*(x)Du(s,x)\big),d^\dagger\hat{B}_s\rangle,\ \
\ \ 0\leq t\leq T.
\end{eqnarray*}

Now let's turn to the existence and uniqueness of solution to the
following infinite horizon BDSDE:
\begin{eqnarray}\label{zz30}
{\rm e}^{-{K'\over2}t}Y_t&=&\int_{t}^{\infty}{\rm e}^{-{K'\over2}s}f(s,Y_s,Z_s)ds+\int_{t}^{\infty}{K'\over2}{\rm e}^{-{K'\over2}s}Y_sds\nonumber\\
&&-\int_{t}^{\infty}{\rm e}^{-{K'\over2}s}\langle
g(s,Y_s,Z_s),d^\dagger\hat{B}_s\rangle-\int_{t}^{\infty}{\rm
e}^{-{K'\over2}s}\langle Z_s,dW_s\rangle,
\end{eqnarray}
or equivalently, for arbitrary $T>0$ and $0\leq t\leq T$,
\begin{eqnarray}
\left\{\begin{array}{l}\label{zz31}
dY_t=-f(t,Y_t,Z_t)dt+\langle g(t,Y_t,Z_t),d^\dagger\hat{B}_t\rangle+\langle Z_t,dW_t\rangle,\nonumber\\
\lim_{T\rightarrow\infty}{\rm e}^{-{K'\over2}T}Y_T=0\ \ \
\rm{a.s.}\nonumber
\end{array}\right.
\end{eqnarray}
We assume that
\begin{description}
\item[(H.1).] Functions $f:\Omega\times[0,\infty)\times\mathbb{R}^1\times\mathbb{R}^{d}{\longrightarrow{\mathbb{R}^1}}$ and $g:\Omega\times[0,\infty)\times\mathbb{R}^1\times\mathbb{R}^{d}{\longrightarrow{\mathbb{R}^{l}}}$ are $\mathscr{F}\otimes\mathscr{B}_{[0,\infty)}\otimes\mathscr{B}_{\mathbb{R}^{1}}\otimes\mathscr{B}_{\mathbb{R}^{d}}$
measurable, and there exist constants $C_1$, $C\geq0$ and
$0\leq\alpha<{1\over2}$ s.t. for any
$(\omega,t)\in\Omega\times[0,\infty)$, $(y_1,z_1)$, $(y_2,z_2)\in
\mathbb{R}^1\times \mathbb{R}^{d}$,
\begin{eqnarray*}
&&|f(t,y_1,z_1)-f(t,y_2,z_2)|^2\leq C_1|y_1-y_2|^2+C |z_1-z_2|^2,\nonumber\\
&&|g(t,y_1,z_1)-g(t,y_2,z_2)|^2\leq C|y_1-y_2|^2+\alpha|z_1-z_2|^2;
\end{eqnarray*}
\item[(H.2).] There exist
constants $K\in\mathbb{R}^+$, $p>d+2$, $K<K'<2K$ and $\mu>0$ with
$2\mu-K'-{p(p+1)\over2}C>0$ s.t. for any
$(\omega,t)\in\Omega\times[0,\infty)$, $y_1$, $y_2\in \mathbb{R}^1$,
$z\in \mathbb{R}^{d}$,
$$(y_1-y_2)(f(t,y_1,z)-f(t,y_2,z))\leq -\mu |y_1-y_2|^2;$$
\item[(H.3).] For $p$, $K$ in (H.2), $f(\cdot, 0, 0)\in
M^{p,-K}([0,\infty); \mathbb{R}^1 )$, $g(\cdot,0,0)\in
M^{p,-K}([0,\infty);\mathbb{R}^{l})$.
\end{description}
\begin{thm}\label{theorem2.3}
Under Conditions {\rm(H.1)}--{\rm(H.3)}, BDSDE (\ref{zz30}) has a
unique solution
\begin{center}
$(Y_{\cdot}, Z_{\cdot})\in S^{p,-K}\bigcap
M^{2,-K}([0,\infty);\mathbb{R}^1)\bigotimes M^{2,-K}([0,\infty);
\mathbb{R}^{d})$,
\end{center}
where the norm in $S^{p,-K}([0,\infty);\mathbb{R}^1)\cap
M^{2,-K}([0,\infty);\mathbb{R}^1)\bigotimes M^{2,-K}([0,\infty);
\mathbb{R}^{d})$ is defined as
\begin{eqnarray*}
\big((E[\sup_{t\geq0}{\rm e}^{-Kt}|\cdot|^p])^{2\over
p}+E[\int_{0}^{\infty}{\rm
e}^{-Kr}|\cdot|^2dr]+E[\int_{0}^{\infty}{\rm
e}^{-Kr}|\cdot|^2dr]\big)^{1\over2},
\end{eqnarray*}
as in Pardoux \cite{pa}.
\end{thm}
{\em Proof}. \underline{Uniqueness}. Let $(Y_{t}^{1},Z_{t}^{1})$ and
$(Y_{t}^{2},Z_{t}^{2})$ be two solutions of BDSDE (4.1). Define
\begin{eqnarray*}
\bar{Y}_t=Y_{t}^{1}-Y_{t}^{2},\ \ \bar{Z}_t=Z_{t}^{1}-Z_{t}^{2},\ \
\ \ \ \ t\geq 0.
\end{eqnarray*}
Applying It$\hat {\rm o}$'s formula to ${\rm
e}^{-Ks}{{|\bar{Y}_s|}^2}$, we have
\begin{eqnarray}\label{zz19}
&&E[{\rm e}^{-Kt}{{|\bar{Y}_t|}^2}]+E[\int_{t}^{T}({1\over
2}-\alpha){\rm e}^{-Ks}|\bar{Z}_s|^2ds]+E[\int_{t}^{T}(2\mu-K-3C){\rm e}^{-Ks}{{|\bar{Y}_s|}^2}ds]\nonumber \\
&\leq&E[{\rm e}^{-KT}{{|\bar{Y}_T|}^2}].
\end{eqnarray}
Taking $K'$ as in Condition (H.2) and noting $2\mu-K'-3C>0$ as well,
we can see that (\ref{zz19}) remains true when $K$ replaced by $K'$.
Therefore, we have
\begin{eqnarray}\label{qi34}
E[{\rm e}^{-K't}{|\bar{Y}_t|}^2]\leq {\rm e}^{-(K'-K)T}E[{\rm
e}^{-KT}{|\bar{Y}_T|}^2].
\end{eqnarray}
Since 
$\sup_{T\geq0}E[{\rm e}^{-KT}{|\bar{Y}_T|}^2]
<\infty$,
taking the limit as $T\to \infty$ in (\ref{qi34}), we have
\begin{eqnarray*}
E[{\rm e}^{-K't}{|\bar{Y}_t|}^2]=0.
\end{eqnarray*}
Then the uniqueness is proved.\\

\underline{Existence}. For each $n\in\mathbb{N}$, we define a
sequence of BDSDEs as follows
\begin{eqnarray}\label{zhao100}
Y_{t}^{n}=\int_{t}^{n}f(s, Y_{s}^{n},
Z_{s}^{n})ds-\int_{t}^{n}\langle g(s,Y_{s}^{n},
Z_{s}^{n}),d^\dagger\hat{B}_s\rangle-\int_{t}^{n}\langle
Z_{s}^{n},dW_s\rangle.
\end{eqnarray}
Let $(Y_{t}^{n}, Z_{t}^{n})_{t\geq n}=(0, 0)$, and according to
Theorem \ref{theorem2.1}, BDSDE (\ref{zhao100}) has a unique
solution $(Y_{\cdot}^{n}, Z_{\cdot}^{n})\in S^{2,-K}\bigcap
M^{2,-K}([0,\infty);\mathbb{R}^1)\bigotimes M^{2,-K}([0,\infty);
\mathbb{R}^{n})$. Also under Conditions (H.1)--(H.3), we can prove
$Y_{\cdot}^{n}\in S^{p,-K}([0,\infty);\mathbb{R}^1)$ in the
following lemma.
\begin{lem}\label{lemma2.4}
Let $(Y_{t}^{n})_{t\geq0}$ be the solution of BDSDE (\ref{zhao100}),
then under Conditions {\rm(H.1)}--{\rm(H.3)}, $Y_{\cdot}^{n}\in
S^{p,-K}([0, \infty);\mathbb{R}^1)$.
\end{lem}
{\em Proof}. Let
\begin{eqnarray*}
&&\psi_M(x)=x^2I_{\{-M\leq x<M\}}+2M(x-M)I_{\{x\geq
M\}}-2M(x+M)I_{\{x<-M\}}\\
&&\varphi_{N,p}(x)=x^{p\over2}I_{\{0\leq
x<N\}}+{p\over2}N^{{p-2}\over2}(x-N)I_{\{x\geq N\}}.
\end{eqnarray*} Applying
generalized It$\hat {\rm o}$'s formula (c.f. Elworthy, Truman and
Zhao \cite{el-tr-zh}) to ${\rm
e}^{-Kr}\varphi_{N,p}\big(\psi_M(Y_r^n)\big)$ to have the following
estimation
\begin{eqnarray}\label{zhang690}
&&{\rm e}^{-Ks}\varphi_{N,p}\big(\psi_M(Y_{s}^{n})\big)-{K}\int_{s}^{n}{\rm e}^{-Kr}\varphi_{N,p}\big(\psi_M(Y_r^n)\big)dr\nonumber\\
&&+{1\over2}\int_{s}^{n}{\rm e}^{-Kr}\varphi^{''}_{N,p}\big(\psi_M(Y_r^n)\big)|\psi_M^{'}(Y_r^n)|^2|Z_r^n|^2dr\nonumber\\
&&+\int_{s}^{n}{\rm e}^{-Kr}\varphi^{'}_{N,p}\big(\psi_M(Y_r^n)\big)I_{\{-M\leq{Y}_r^n<M\}}|{Z}_r^n|^2dr\nonumber\\
&\leq&\int_{s}^{n}{\rm e}^{-Kr}\varphi^{'}_{N,p}\big(\psi_M(Y_r^n)\big)\psi_M^{'}(Y_r^n){f}(r,Y_r^n,Z_r^n)dr\nonumber\\
&&+\int_{s}^{n}{\rm e}^{-Kr}\varphi^{'}_{N,p}\big(\psi_M(Y_r^n)\big)I_{\{-M\leq{Y}_r^n<M\}}|g(r,Y_r^n,Z_r^n)|^2dr\nonumber\\
&&+{1\over2}\int_{s}^{n}{\rm e}^{-Kr}\varphi^{''}_{N,p}\big(\psi_M(Y_r^n)\big)|\psi_M^{'}(Y_r^n)|^2|g(r,Y_r^n,Z_r^n)|^2dr\nonumber\\
&&-\int_{s}^{n}\langle{\rm
e}^{-Kr}\varphi^{'}_{N,p}\big(\psi_M(Y_r^n)\big)\psi_M^{'}(Y_r^n)g(r,Y_r^n,
Z_r^n),d^\dagger\hat{B}_r\rangle\nonumber\\
&&-\int_{s}^{n}\langle{\rm
e}^{-Kr}\varphi^{'}_{N,p}\big(\psi_M(Y_r^n)\big)\psi_M^{'}(Y_r^n){Z}_r^n,dW_r\rangle.
\end{eqnarray}
As $(Y_{\cdot}^{t,\cdot},Z_{\cdot}^{t,\cdot})\in S^{2,-K}\bigcap
M^{2,-K}([0,\infty);{\mathbb{R}^{1}})\bigotimes
M^{2,-K}([0,\infty);{\mathbb{R}^{d}})$ and
$\varphi^{'}_{N,p}\big(\psi_M(Y_r^n)\big)\psi_M^{'}(Y_r^n)$ is
bounded, taking the expectation on both sides, we know that all the
stochastic integrals have zero expectation. Using Conditions
(H.1)-(H.3) and taking first the limit as $M\to \infty$, then the
limit as $N\to \infty$, by the monotone convergence theorem, we have
\begin{eqnarray}\label{zhang2}
&&\big(p\mu-{K}-{{p(p+1)}\over2}C-(3+{{p(p-1)}\over2}C)\varepsilon\big)E[\int_{s}^{\infty}{\rm e}^{-Kr}{|{Y}_r^n|}^pr]\nonumber\\
&&+{p\over4}\big(2p-3-(2p-2)\alpha-(2p-2)\alpha\varepsilon\big)E[\int_{s}^{\infty}{\rm e}^{-Kr}{{|{Y}_r^n|}^{p-2}}|{Z}_r^n|^2dr]\nonumber\\
&\leq&C_pE[\int_{0}^{\infty}{\rm
e}^{-Kr}|f(r,0,0)|^pdr]+C_pE[\int_{0}^{\infty}{\rm
e}^{-Kr}|g(r,0,0)|^pdr]<\infty.
\end{eqnarray}
Note that here and in the following the constant $\varepsilon$ can
be chosen to be sufficiently small and $C_p$ is a generic constant.
Due to Conditions (H.1), (H.2) and the arbitrariness of
$\varepsilon$, all the terms on the left hand side of (\ref{zhang2})
are positive. Furthermore, by the B-D-G inequality, Cauchy-Schwartz
inequality and Young inequality, from (\ref{zhang690}) we have
\begin{eqnarray}\label{zhang692}
&&E[\sup_{s\geq0}{\rm e}^{-Ks}{{|{Y}_s^{n}|}^p}]\nonumber\\
&\leq&C_pE[\int_{0}^{\infty}{\rm
e}^{-Kr}{{|{Y}_r^n|}^{p-2}}|{Z}_r^n|^2dr]+C_pE[\int_{0}^{\infty}{\rm
e}^{-Kr}|Y_r^n|^pdr]\nonumber\\
&&+C_pE[\sqrt{\int_{s}^{\infty}\big({\rm
e}^{-Kr}\varphi^{'}_{N,p}\big(\psi_M(Y_r^n)\big)|\psi_M^{'}(Y_r^n)|^2\big)\big({\rm
e}^{-Kr}\varphi^{'}_{N,p}\big(\psi_M(Y_r^n)\big)|g(r,Y_r^n,
Z_r^n)|^2\big)dr}]\nonumber\\
&&+C_pE[\sqrt{\int_{s}^{\infty}\big({\rm
e}^{-Kr}\varphi^{'}_{N,p}\big(\psi_M(Y_r^n)\big)|\psi_M^{'}(Y_r^n)|^2\big)\big({\rm
e}^{-Kr}\varphi^{'}_{N,p}\big(\psi_M(Y_r^n)\big)|{Z}_r^n|^2\big)dr}]\nonumber\\
&\leq&C_pE[\int_{0}^{\infty}{\rm
e}^{-Kr}{{|{Y}_r^n|}^{p-2}}|{Z}_r^n|^2dr]+C_pE[\int_{0}^{\infty}{\rm
e}^{-Kr}|Y_r^n|^pdr]\nonumber\\
&&+\varepsilon E[\sup_{s\geq0}\big({\rm
e}^{-Ks}\varphi^{'}_{N,p}\big(\psi_M(Y_s^n)\big)|\psi_M^{'}(Y_s^n)|^2\big)]+C_pE[\int_{0}^{\infty}{\rm
e}^{-Kr}\varphi^{'}_{N,p}\big(\psi_M(Y_r^n)\big)|g(r,Y_r^n,
Z_r^n)|^2dr]\nonumber\\
&&+C_pE[\int_{0}^{\infty}{\rm
e}^{-Kr}\varphi^{'}_{N,p}\big(\psi_M(Y_r^n)\big)|Z_r^n|^2dr].
\end{eqnarray}
Taking the limits as $M$, $N\rightarrow\infty$ and applying the
monotone convergence theorem, we have
\begin{eqnarray}\label{zhang4}
E[\sup_{s\geq0}{\rm
e}^{-Kt}|Y_{s}^{n}|^p]&\leq&C_pE[\int_{0}^{\infty}{\rm
e}^{-Kr}|Y_r^n|^{p-2}|Z_r^n|^2dr]+C_pE[\int_{0}^{\infty}{\rm
e}^{-Kr}|Y_r^n|^pdr].
\end{eqnarray}
By (\ref{zhang2}), $Y_{\cdot}^{n}\in
S^{p,-K}([0,\infty);\mathbb{R}^1)$. Lemma \ref{lemma2.4} is proved.
$\hfill\diamond$
\\

\begin{rmk}\label{remark2.5}
The proof of Lemma \ref{lemma2.4} also works with $p$ replaced by
$2$. Note that if $f(\cdot,0,0)\in
M^{p,-K}([0,\infty);\mathbb{R}^1)$, then by H$\ddot{\textrm{o}}$lder
inequality, it turns out that
$f(\cdot,0,0)\in M^{2,-K}([0,\infty);\mathbb{R}^1)$ and
$g(\cdot,0,0)\in M^{2,-K}([0,\infty);\mathbb{R}^l)$. So it is easy
to see in (\ref{zhang2}) with $p$ replaced by $2$ that
\begin{center}
${(Y_{\cdot}^{n},Z_{\cdot}^{n})}\in
M^{2,-K}([0,\infty);\mathbb{R}^1)\bigotimes
M^{2,-K}([0,\infty);\mathbb{R}^{d})$.
\end{center}
For the rest of our paper, we will leave out the similar
localization argument as in the proof of Lemma \ref{lemma2.4} when
applying It$\hat {\rm o}$'s formula to save the space of this paper.
\end{rmk}

Then back to the proof of Theorem \ref{theorem2.3}. We will show
that ${(Y_{\cdot}^{n}, Z_{\cdot}^{n})}$ is a Cauchy sequence in the
space of $S^{p,-K}([0,\infty);\mathbb{R}^1)\cap
M^{2,-K}([0,\infty);\mathbb{R}^1)\bigotimes M^{2,-K}([0,\infty);
\mathbb{R}^{d})$. First we show that, for $m,n\in\mathbb{N}$ and
$m\geq n$,
\begin{eqnarray*}
\lim_{n,m\rightarrow\infty}E[\sup_{t\geq0}{\rm
e}^{-Kt}|Y_{t}^{m}-Y_{t}^{n}|^p]=0.
\end{eqnarray*}
Define $\bar{Y}_{t}^{m,n}={Y}_{t}^{m}-{Y}_{t}^{n}$,
$\bar{Z}_{t}^{m,n}={Z}_{t}^{m}-{Z}_{t}^{n}$.

($\textrm{i}$) When $n\leq t\leq m$,
\begin{eqnarray*}
\bar{Y}_{t}^{m,n}={Y}_{t}^{m}=\int_{t}^{m}f(s,Y_{s}^{m},Z_{s}^{m})ds-\int_{t}^{m}\langle
g(s,Y_{s}^{m},Z_{s}^{m}),d^\dagger\hat{B}_s\rangle-\int_{t}^{m}\langle
Z_{s}^{m},dW_s\rangle.
\end{eqnarray*}
Some similar calculations as in (\ref{zhang2}) and (\ref{zhang4})
lead to
\begin{eqnarray}\label{zhang6}
E[\sup_{n\leq t\leq m}{\rm e}^{-Kt}|Y_{t}^{m}|^p]&\leq&C_pE[\int_{n}^{m}{\rm e}^{-Kr}|Y_{r}^{m}|^{p-2}|Z_{r}^{m}|^2dr]+C_pE[\int_{n}^{m}{\rm e}^{-Kr}|Y_{r}^{m}|^pdr]\\
&&+C_pE[\int_{n}^{m}{\rm e}^{-Kr}(|f(r,0,0)|^p+|g(r,0,0)|^p)dr]
\longrightarrow0,\ {\rm as}\ n,\ m\longrightarrow\infty.\nonumber
\end{eqnarray}

$(\textrm{i}$$\textrm{i})$ When $0\leq t\leq n$,
\begin{eqnarray*}
\bar{Y}_{t}^{m,n}&=&{Y}_{n}^{m}+\int_{t}^{n}f(r,Y_{r}^{m},Z_{r}^{m})-f(r,Y_r^n,Z_r^n)\\
&&-\int_{t}^{n}\langle
g(r,Y_{r}^{m},Z_{r}^{m})-g(r,Y_r^n,Z_r^n),d^\dagger\hat{B}_r\rangle-\int_{t}^{n}\langle\bar{Z}_{r}^{m,n},dW_r\rangle.
\end{eqnarray*}
Applying It$\hat {\rm o}$'s formula to ${\rm
e}^{-Kr}|\bar{Y}_{r}^{m,n}|^p$ and following a similar calculation
as in (\ref{zhang690}) and (\ref{zhang2}), we have for $s\leq n$,
\begin{eqnarray}\label{zhang8}
E[\int_{0}^{n}{\rm
e}^{-Kr}|\bar{Y}_{r}^{m,n}|^{p-2}|\bar{Z}_{r}^{m,n}|^2dr]+E[\int_{0}^{n}{\rm
e}^{-Kr}|\bar{Y}_{r}^{m,n}|^pdr]\leq C_pE[{\rm
e}^{-Kn}|Y_{n}^{m}|^p].
\end{eqnarray}
From $(\textrm{i})$, the right hand side of the above inequality
converges to $0$ as $n$, $m\longrightarrow\infty$.
By some similar calculations as in (\ref{zhang4}), we have
\begin{eqnarray*}
E[\sup_{0\leq t\leq n}{\rm e}^{-Kt}|\bar{Y}_{t}^{m,n}|^p]\leq
C_pE[{\rm e}^{-Kn}|Y_{n}^{m}|^p]\longrightarrow0\ \ \ {\rm as} \
n,m\longrightarrow\infty.
\end{eqnarray*}
From $(\textrm{i})$ $(\textrm{i}$$\textrm{i})$, we have for
$m,n\in\mathbb{N}$,
\begin{eqnarray*}
\lim_{n,m\rightarrow\infty}E[\sup_{t\geq0}{\rm
e}^{-Kt}|Y_{t}^{m}-Y_{t}^{n}|^p]=0.
\end{eqnarray*}
It is easy to see that the above arguments also hold for $p=2$ in
(\ref{zhang6}) and (\ref{zhang8}). Noting Remark \ref{remark2.5}, we
have as $n$, $m\longrightarrow\infty$
\begin{eqnarray*}
E[\int_{0}^{\infty}{\rm
e}^{-Kr}|\bar{Y}_{r}^{m,n}|^2dr]+E[\int_{0}^{\infty}{\rm
e}^{-Kr}|\bar{Z}_{r}^{m,n}|^2dr]\longrightarrow0.
\end{eqnarray*}
Therefore, ${(Y_{\cdot}^{n}, Z_{\cdot}^{n})}$ is a Cauchy sequence in
the Banach space $S^{p,-K}([0,\infty);\mathbb{R}^1)\cap
M^{2,-K}([0,\infty);\mathbb{R}^1)\bigotimes\\
M^{2,-K}([0,\infty);\mathbb{R}^{d})$.

We take $(Y_t, Z_t)_{t\geq0}$ as the limit of
$(Y_{t}^{n},Z_{t}^{n})_{t\geq0}$ in
$S^{p,-K}([0,\infty);\mathbb{R}^1)\cap
M^{2,-K}([0,\infty);\mathbb{R}^1)\bigotimes
M^{2,-K}\\([0,\infty);\mathbb{R}^{d})$ and then show that $(Y_t,
Z_t)_{t\geq0}$ is the solution of BDSDE (\ref{zz30}). First note that
for $t\leq n$, (\ref{zhao100}) is equivalent to
\begin{eqnarray}\label{zhang9}
{\rm e}^{-{K'\over2}t}Y_{t}^{n}&=&\int_{t}^{n}{\rm e}^{-{K'\over2}s}f(s,Y_{s}^{n},Z_{s}^{n})ds+\int_{t}^{n}{K'\over2}{\rm e}^{-{K'\over2}s}Y_{s}^{n}ds\nonumber\\
&&-\int_{t}^{n}{\rm e}^{-{K'\over2}s}\langle
g(s,Y_{s}^{n},Z_{s}^{n}),d^\dagger\hat{B}_s\rangle-\int_{t}^{n}{\rm
e}^{-{K'\over2}s}\langle Z_{s}^{n},dW_s\rangle.
\end{eqnarray}
Actually BDSDE (\ref{zhang9}) converges to BDSDE (\ref{zz30}) in
$L^2(\Omega)$ as $n\longrightarrow\infty$. To see this, we verify
the convergence term by term. For the first term,
\begin{eqnarray*}
&&E[\ |{\rm e}^{-{K'\over2}t}Y_{t}^{n}-{\rm
e}^{-{K'\over2}t}Y_t|^2]\leq E[\sup_{t\geq0}{\rm
e}^{-Kt}|Y_{t}^{n}-Y_t|^2]\longrightarrow0.
\end{eqnarray*}
For the second term, by H$\ddot{\textrm{o}}$lder inequality,
\begin{eqnarray*}
&&E[\ |\int_{t}^{n}{\rm e}^{-{K'\over2}s}f(s,Y_{s}^{n},Z_{s}^{n})ds-\int_{t}^{\infty}{\rm e}^{-{K'\over2}s}f(s,Y_s,Z_s)ds|^2]\\
&\leq&2E[\int_{t}^{n}{\rm e}^{-(K'-K)s}ds\int_{t}^{n}{\rm e}^{-Ks}|f(s,Y_{s}^{n},Z_{s}^{n})-f(s,Y_s,Z_s)|^2ds]\\
&&+2E[\int_{n}^{\infty}{\rm e}^{-(K'-K)s}ds\int_{n}^{\infty}{\rm
e}^{-Ks}|f(s,Y_s,Z_s)|^2ds]\longrightarrow0.
\end{eqnarray*}
We can deal with the third term similarly as above and deal with two
stochastic integration terms by It$\hat {\rm o}$'s isometry.
Thus ${(Y_t, Z_t)}_{t\geq 0}$ is the solution of BDSDE (\ref{zz30})
and the proof of Theorem \ref{theorem2.3} is completed.
$\hfill\diamond$\\

Then we consider the existence and uniqueness of solution to the
following infinite horizon forward BDSDE:
\begin{eqnarray}\label{zz32}
{\rm e}^{-{K'\over2}s}Y_s^{t,x}&=&\int_{s}^{\infty}{\rm e}^{-{K'\over2}r}f(X_r^{t,x},Y_r^{t,x},Z_r^{t,x})dr+\int_{s}^{\infty}{K'\over2}{\rm e}^{-{K'\over2}r}Y_r^{t,x}dr\\
&&-\int_{s}^{\infty}{\rm e}^{-{K'\over2}r}\langle
g(X_r^{t,x},Y_r^{t,x},Z_r^{t,x}),d^\dagger\hat{B}_r\rangle-\int_{s}^{\infty}{\rm
e}^{-{K'\over2}r}\langle Z_r^{t,x},dW_r\rangle,\ \ s\geq0.\nonumber
\end{eqnarray}
We replace Condition (A.1) by
\begin{description}
\item[(A.1)$^*$.] Functions $f: \mathbb{R}^d\times\mathbb{R}^1\times\mathbb{R}^{d}{\longrightarrow{\mathbb{R}^1}}$ and $g: \mathbb{R}^d\times\mathbb{R}^1\times\mathbb{R}^{d}{\longrightarrow{\mathbb{R}^{l}}}$ are $\mathscr{B}_{\mathbb{R}^{d}}\otimes\mathscr{B}_{\mathbb{R}^{1}}\otimes\mathscr{B}_{\mathbb{R}^{d}}$
measurable, and there exist constants $C_0$, $C_1$, $C\geq0$ and
$0\leq\alpha<{1\over2}$ s.t. for any $(x_1, y_1, z_1)$, $(x_2, y_2,
z_2)\in \mathbb{R}^d\times \mathbb{R}^1\times \mathbb{R}^{d}$,
\begin{eqnarray*}
&&|f(x_1, y_1, z_1)-f(x_2, y_2, z_2)|^2\leq C_0|x_1-x_2|^2+C_1|y_1-y_2|^2+C|z_1-z_2|^2,\nonumber\\
&&|g(x_1, y_1, z_1)-g(x_2, y_2, z_2)|^2\leq
C_0|x_1-x_2|^2+C|y_1-y_2|^2+\alpha|z_1-z_2|^2.
\end{eqnarray*}
\end{description}
\begin{prop}\label{zz34}
Under Conditions {\rm(A.1)$^*$}, {\rm(A.3)}, {\rm(A.4)}, BDSDE
(\ref{zz32}) has a unique solution
\begin{center}
$(Y_{\cdot}^{t,x}, Z_{\cdot}^{t,x})\in
S^{p,-K}([0,\infty);\mathbb{R}^1) \cap
M^{2,-K}([0,\infty);\mathbb{R}^1)\bigotimes
M^{2,-K}([0,\infty);\mathbb{R}^{d})$.
\end{center}
\end{prop}
\begin{rmk}\label{qi315} 
For $s\in[0,t]$, BDSDE (\ref{zz32}) is equivalent to the following
BDSDE
\begin{eqnarray*}\label{zhang663}
{\rm e}^{-{K'\over2}s}Y_s^{x}&=&\int_{s}^{\infty}{\rm e}^{-{K'\over2}r}f(x,Y_r^{x},Z_r^{x})dr+\int_{s}^{\infty}{K'\over2}{\rm e}^{-{K'\over2}r}Y_r^{x}dr\nonumber\\
&&-\int_{s}^{\infty}{\rm e}^{-{K'\over2}r}\langle
g(x,Y_r^{x},Z_r^{x}),d^\dagger\hat{B}_r\rangle-\int_{s}^{\infty}{\rm
e}^{-{K'\over2}r}\langle Z_r^{x},dW_r\rangle.
\end{eqnarray*}
To unify the notation, we define
$({Y}_s^{t,x},{Z}_s^{t,x})=({Y}_s^{x},{Z}_s^{x})$ when $s\in[0,t)$.
\end{rmk}
{\em Proof of Proposition \ref{zz34}}. Let
\begin{eqnarray*} \hat{f}(s,y,z)=f(X_{s}^{t,x},y,z),\ \ \
\hat{g}(s,y,z)=g(X_{s}^{t,x},y,z).
\end{eqnarray*}
We need to verify that $\hat{f}$, $\hat{g}$ satisfy Conditions
(H.1)--(H.3) in Theorem \ref{theorem2.3}. It is obvious that
$\hat{f}$, $\hat{g}$ satisfy (H.1) and (H.2), so we only need to
show that $\hat{f}$, $\hat{g}$ satisfy (H.3) as well, i.e.
\begin{eqnarray*}
E[\int_{0}^{\infty}{\rm e}^{-Ks}|\hat{f}(s,0,0)|^pds]<\infty\ {\rm
and}\ E[\int_{0}^{\infty}{\rm e}^{-Ks}|\hat{g}(s,0,0)|^pds]<\infty.
\end{eqnarray*}
Since
\begin{eqnarray*}
E[\int_{0}^{\infty}{\rm e}^{-Ks}|\hat{f}(s,0,0)|^pds]
\leq C_pE[\int_{0}^{\infty}{\rm
e}^{-Ks}C_0^p|X_{s}^{t,x}|^pds]+C_pE[\int_{0}^{\infty}{\rm
e}^{-Ks}|f(0,0,0)|^pds],
\end{eqnarray*}
we only need to prove $E[\int_{0}^{\infty}{\rm
e}^{-Ks}|X_{s}^{t,x}|^pds]<\infty$. Now applying It$\hat {\rm o}$'s
formula to ${\rm e}^{-Kr}|X_{r}^{t,x}|^p$ and noticing Condition
(A.4), we have
\begin{eqnarray*}
E[\int_{t}^{s}{\rm e}^{-Kr}|X_{r}^{t,x}|^pdr]\leq {\rm
e}^{-Kt}|x|^p+C_pE[\int_{t}^{s}{\rm
e}^{-Kr}(|b(0)|^p+\|\sigma(0)\|^p)dr]<\infty.
\end{eqnarray*}
Taking the limit of $s$ and noting that $(X_{s}^{t,x})_{s<t}=x$, we
have
$E[\int_{0}^{\infty}{\rm e}^{-Kr}|X_{r}^{t,x}|^pdr]<\infty$.
So $E[\int_{0}^{\infty}{\rm e}^{-Ks}|\hat{f}(s,0,0)|^pds]<\infty$.
Similarly, $E[\int_{0}^{\infty}{\rm
e}^{-Ks}|\hat{g}(s,0,0)|^pds]<\infty$. $\hfill\diamond$\\

Now we prove the other assumption in Theorem \ref{zz42}, i.e. the
regularity of solutions of infinite horizon BDSDEs. An simple
application of stochastic flow property proved in \cite{ku2} leads
to
\begin{lem}\label{zz35} Under Condition {\rm(A.4)}, for arbitrary $T$ and $t$, $t'\in[0,T]$,
$x$, $x'$ belonging to an arbitrary bounded set in $\mathbb{R}^d$,
the diffusion process $(X_{s}^{t,x})_{s\geq0}$ defined in SDE
(\ref{qi17}) satisfies
\begin{eqnarray*}
E[\int_{0}^{\infty}{\rm
e}^{-Kr}|X_r^{t',x'}-X_r^{t,x}|^pdr]&\leq&C_p(|x'-x|^p+|t'-t|^{p\over2})\
\ \ \rm{a.s.}
\end{eqnarray*}
\end{lem}
$\hfill\diamond$\\

We concentrate ourselves on the regularity of infinite horizon BDSDE
(\ref{zz36}), which is
a simpler form of BDSDE (\ref{zz32}). For arbitrary given terminal
time $T$, the form of BDSDE (\ref{zz36}) on $[t,T]$ is (\ref{zz22}).
\begin{prop}\label{zz39} Under Conditions {\rm(A.1)}--{\rm(A.4)}, let $(Y_{s}^{t,x})_{s\geq0}$ be the solution of BDSDE (\ref{zz36}), then for arbitrary $T$ and $t\in[0,T]$, $x\in\mathbb{R}^d$, $(t,x)\longrightarrow Y_{t}^{t,x}$ is a.s.
continuous.
\end{prop}
{\em Proof}. For $t$, $t'$, $r\geq0$, let
\begin{eqnarray*}
&&\bar{Y}_r=Y_{r}^{t',x'}-Y_{r}^{t,x},\ \ \
\bar{Z}_r=Z_{r}^{t',x'}-Z_{r}^{t,x}.
\end{eqnarray*}
Applying It$\hat {\rm o}$'s formula to ${\rm
e}^{-{{pK'}\over2}r}|\bar{Y}_r|^p$ and following a similar
calculation as in (\ref{zhang690}), we have
for $0\leq s\leq T$,
\begin{eqnarray}\label{zhang13}
&&{\rm e}^{-{{pK'}\over2}s}|\bar{Y}_s|^p+(p\mu-{{pK'}\over2}-{{p(p+1)}\over2}C-\varepsilon)\int_{s}^{T}{\rm e}^{-{{pK'}\over2}r}|\bar{Y}_{r}|^pdr\nonumber\\
&&+{{p(2p-3)}\over4}\int_{s}^{T}{\rm e}^{-{{pK'}\over2}r}|\bar{Y}_{r}|^{p-2}|\bar{Z}_r|^2dr\nonumber\\
&\leq&{\rm e}^{-{{pK'}\over2}T}|\bar{Y}_T|^p+C_p\int_{s}^{T}{\rm e}^{-{{pK'}\over2}r}|\bar{X}_{r}|^pdr-p\int_{s}^{T}{\rm e}^{-{{pK'}\over2}r}|\bar{Y}_{r}|^{p-2}\bar{Y}_{r}\langle\bar{g}_r,d^\dagger\hat{B}_r\rangle\nonumber\\
&&-p\int_{s}^{T}{\rm
e}^{-{{pK'}\over2}r}|\bar{Y}_{r}|^{p-2}\bar{Y}_{r}\langle\bar{Z}_r,dW_r\rangle.
\end{eqnarray}
Noticing Condition (A.3), 
for $0\leq s\leq T$, we have
\begin{eqnarray}\label{zhang14}
&&E[{\rm e}^{-{{pK'}\over2}s}|\bar{Y}_s|^p]+E[\int_{s}^{T}{\rm e}^{-{{pK'}\over2}r}|\bar{Y}_{r}|^{p}dr]+E[\int_{s}^{T}{\rm e}^{-{{pK'}\over2}r}|\bar{Y}_{r}|^{p-2}|\bar{Z}_r|^2dr]\nonumber\\
&\leq&C_pE[{\rm
e}^{-{{pK'}\over2}T}|\bar{Y}_T|^p]+C_pE[\int_{s}^{T}{\rm
e}^{-{{pK'}\over2}r}|\bar{X}_{r}|^pdr].
\end{eqnarray}
Since 
$E[{\rm e}^{-{{pK'}\over2}T}|\bar{Y}_T|^p]\leq E[\sup_{s\geq0}{\rm
e}^{-Ks}|\bar{Y}_s|^p]<\infty$,
by the Lebesgue's dominated convergence theorem, we have
\begin{eqnarray}\label{zhang15}
\lim_{T\rightarrow\infty}E[{\rm
e}^{{{pK'}\over2}T}|\bar{Y}_T|^p]=E[(\lim_{T\rightarrow\infty}{\rm
e}^{-{K'\over2}T}|\bar{Y}_T|)^p]=0.
\end{eqnarray}
So taking the limit of $T$ in (\ref{zhang14}), by Lemma \ref{zz35}
and the monotone convergence theorem, we have
\begin{eqnarray}\label{zhang16}
E[\int_{0}^{\infty}{\rm
e}^{-{{pK'}\over2}r}|\bar{Y}_{r}|^{p-2}|\bar{Z}_r|^2dr]+E[\int_{0}^{\infty}{\rm
e}^{-{{pK'}\over2}r}|\bar{Y}_{r}|^pdr]\leq
C_pE[\int_{0}^{\infty}{\rm
e}^{-Kr}|\bar{X}_{r}|^pdr].
\end{eqnarray}
From (\ref{zhang13}), by B-D-G inequality and (\ref{zhang15}), we
have
\begin{eqnarray*}
&&E[\sup_{s\geq0}{\rm e}^{-{{pK'}\over2}s}|\bar{Y}_s|^p]\\
&\leq&C_pE[\int_{0}^{\infty}{\rm
e}^{-{{pK'}\over2}r}|\bar{X}_{r}|^pdr]+C_pE[\int_{0}^{\infty}{\rm
e}^{-{{pK'}\over2}r}|\bar{Y}_{r}|^pdr]+C_pE[\int_{0}^{\infty}{\rm
e}^{-{{pK'}\over2}r}|\bar{Y}_{r}|^{p-2}|\bar{Z}_r|^2dr].
\end{eqnarray*}
By the above inequality, Lemma \ref{zz35} and (\ref{zhang16}), for
arbitrary $T>0$, $t$, $t'\in[0,T]$, $x$, $x'$ belonging to an
arbitrary bounded set in $\mathbb{R}^d$, we have
\begin{eqnarray}\label{zhang18}
E[\sup_{s\geq0}{\rm e}^{-{pK}s}|\bar{Y}_s|^p]\leq
C_pE[\int_{s}^{T}{\rm e}^{-{{pK'}\over2}r}|\bar{X}_{r}|^pdr]\leq
C_p(|x'-x|^p+|t'-t|^{p\over2}).
\end{eqnarray}
Noting $p>d+2$ in (\ref{zhang18}), by Kolmogorov Lemma (see e.g.
\cite{ku2}), we have $Y_{s}^{(\cdot,\cdot)}$ has a continuous
modification for $t\in[0,T]$ and $x$ belonging to an arbitrary
bounded set in $\mathbb{R}^d$ under the norm $\sup_{s\geq0}{\rm
e}^{-{K}s}|Y_{s}^{(\cdot,\cdot)}|$. In particular,
\begin{eqnarray*}
\lim_{t'\rightarrow t\atop x'\rightarrow x}{\rm
e}^{-{K}t'}|Y_{t'}^{t',x'}-Y_{t'}^{t,x}|=0.
\end{eqnarray*}
Thus we have a.s.
\begin{eqnarray*}
\lim_{t'\rightarrow t\atop x'\rightarrow x}|{\rm
e}^{-{K}t'}Y_{t'}^{t',x'}-{\rm
e}^{-{K}t}Y_{t}^{t,x}|\leq\lim_{t'\rightarrow t\atop x'\rightarrow
x}(|{\rm e}^{-{K}t'}Y_{t'}^{t',x'}-{\rm
e}^{-{K}t'}Y_{t'}^{t,x}|+|{\rm e}^{-{K}t'}Y_{t'}^{t,x}-{\rm
e}^{-{K}t}Y_{t}^{t,x}|)=0.
\end{eqnarray*}
The convergence of the second term follows from the continuity of
$Y_{s}^{t,x}$ in $s$. That is to say ${\rm e}^{-{K}t}Y_{t}^{t,x}$ is
a.s. continuous, therefore $Y_{t}^{t,x}$ is continuous w.r.t.
$t\in[0,T]$ and $x$ belonging to an arbitrary bounded set in
$\mathbb{R}^d$.

Denote by $\bar{B}(0,R)$ the closed ball in $\mathbb{R}^d$ of radius
$R$ centered at $0$. It is obvious that
$\bigcup_{R=1}^{\infty}\bar{B}(0,R)=\mathbb{R}^d$. $Y_{t}^{t,x}$ is
continuous w.r.t $t\in[0,T]$ and $x\in\bar{B}(0,R)$ on $\Omega^R$.
Take $\tilde{\Omega}=\bigcap_{R=1}^{\infty}{\Omega}^{R}$, then
$P(\tilde{\Omega})=1$. Now for any $t\in[0,T]$ and
$x\in\mathbb{R}^d$, there exists an $R$ s.t. $x\in\bar{B}(0,R)$. On
the other hand, for all $\omega\in\tilde{\Omega}$, it is obvious
that $\omega\in\Omega^{R}$. So $Y_{t}^{t,x}$ is continuous w.r.t.
$t\in[0,T]$ and $x\in\mathbb{R}^{d}$ on $\tilde{\Omega}$.
Proposition \ref{zz39} is proved. $\hfill\diamond$\\

\section{Stationary Property of Stochastic Viscosity Solutions of SPDEs}\label{s27}
\setcounter{equation}{0}

\ \ \ \ With the regularity of solution of BDSDE (\ref{zz36}), for
arbitrary given $T$, we can obtain a stochastic viscosity solution
of SPDE (\ref{zz20}) on the time interval $[0,T]$ through BDSDE
(\ref{zz36}).
\begin{thm}\label{zz37}
Under Conditions {\rm(A.1)}--{\rm(A.4)}, for arbitrary given $T$ and
$t\in[0,T]$, $x\in\mathbb{R}^d$, let $v(t,x)\triangleq
Y_{T-t}^{T-t,x}$, where $(Y_{s}^{t,x},Z_{s}^{t,x})$ is the solution
of BDSDE (\ref{zz36}) with $\hat{B}_s={B}_{T-s}-{B}_T$ for all
$s\geq0$. Then $v(t,x)$ is continuous w.r.t. $t$ and $x$ and is a
stochastic viscosity solution of SPDE (\ref{zz20}) on the time
interval $[0,T]$.
\end{thm}
{\em Proof}. Notice that Condition (A.1) is stronger than (A.2)$'$,
so by Theorem \ref{zz38} BDSDE (\ref{zz36}) has a unique solution
$(Y_{\cdot}^{t,x}, Z_{\cdot}^{t,x})\in
S^{p,-K}([0,\infty);\mathbb{R}^1) \cap
M^{2,-K}([0,\infty);\mathbb{R}^1)\bigotimes
M^{2,-K}([0,\infty);\mathbb{R}^{d})$. On $[t,T]$, BDSDE (\ref{zz36})
has a form of (\ref{zz22}) which can be associated with SPDE
(\ref{zz20}) on $[0,T]$ through time reversal transformation in
(\ref{zz21}). First note that by Proposition \ref{zz39}, $v(t,x)$
defined by $Y_{T-t}^{T-t,x}$ is a.s. continuous w.r.t. $t\in[0,T]$
and $x\in\mathbb{R}^d$. Moreover, since
$X_s^{T,X_T^{t,x}}=X_s^{t,x}$ for $s\geq T$, by the uniqueness of
BDSDE (\ref{zz36}) we have
$Y_T^{T,X_T^{t,x}}=Y_T^{t,x}$ a.s., where $Y_\cdot^{T,x}
$ is the solution of BDSDE (\ref{zz36}) when the diffusion process
$X$ defined in (\ref{qi17})
starts at time $T$ and point $x\in\mathbb{R}^{d}$. Therefore 
$E[|v(0,X_T^{t,x})|^2]=E[|Y_T^{T,X_T^{t,x}}|^2]=E[|Y_T^{t,x}|^2]<\infty$.
By Theorem \ref{theorem2.2} and Remark \ref{zz24}, we know that
$v(t,x)$ is a stochastic viscosity solution of SPDE (\ref{zz20}) on
the time interval $[0,T]$. Theorem \ref{zz37} is proved.
$\hfill\diamond$
\\

In the following, we show that the $v(t,x)$ constructed in Theorem
\ref{zz37} is a stationary solution of SPDE (\ref{zz20}). For this,
we need first prove a claim that
$v(t,x)(\omega)=Y_{T-t}^{T-t,x}(\hat{\omega})$ is independent of the
choice of $T$. This independence can be proved by a similar argument
as in \cite{zh-zh1} (Page 186-187) since it is unrelated to which
kind of solution (weak solution or stochastic viscosity solution)
$v$ is. Therefore, for any $T'\geq T$,
$Y_{T-t}^{T-t,x}(\hat{\omega})=Y_{T'-t}^{T'-t,x}({\hat{\omega}}')$
when $0\leq t\leq T$, where $\hat{\omega}(s)={B}_{T-s}-{B}_{T}$ and
${\hat{\omega}}'(s)={B}_{T'-s}-{B}_{T'}$.

On the probability space $(\Omega,\mathscr{F},P)$, we define
${\theta}_{t}=(\hat{\theta}_{t})^{-1}$, $t\geq0$. Actually $\hat{B}$
is a two-sided Brownian motion, so
$(\hat{\theta}_{t})^{-1}=\hat{\theta}_{-t}$ is well defined (see
\cite{ar}). It is easy to see that ${\theta}_{t}$ is a shift w.r.t.
${B}$ satisfying
\begin{description}
\item[$(\textrm{i})$]$P\cdot({\theta}_{t})^{-1}=P$;
\item[$(\textrm{i}\textrm{i})$]${\theta}_{0}=I$;
\item[$(\textrm{i}\textrm{i}\textrm{i})$]${\theta}_{s}\circ{\theta}_{t}={\theta}_{s+t}$;
\item[$(\textrm{iv})$]${\theta}_{t}\circ{B}_s={B}_{s+t}-{B}_{t}$.
\end{description}
By Theorem \ref{zz42} and the relationship between $\theta$ and
$\hat{\theta}$, we have
\begin{eqnarray*}
{\theta}_rv(t,x)(\omega)=\hat{\theta}_{-r}Y_{T-t}^{T-t,x}(\hat{\omega})=\hat{\theta}_{-r}\hat{\theta}_{r}Y_{T-t-r}^{T-t-r,x}(\hat{\omega})=Y_{T-t-r}^{T-t-r,x}(\hat{\omega})=v(t+r,x)(\omega),
\end{eqnarray*}
for all $r\geq0$ and $T\geq t+r$, $x\in\mathbb{R}^{d}$ a.s. In
particular, let
$Y(x,\omega)=v(0,x)(\omega)=Y_{T}^{T,x}(\hat{\omega})$, then the
above formula implies (\ref{zhao001}):
\begin{eqnarray*}
{\theta}_tY(x,\omega)=Y(x,{\theta}_t\omega)=v(t,x)(\omega)=v(t,x,v(0,x)(\omega))(\omega)=v(t,x,Y(x,\omega))(\omega),\
\end{eqnarray*}
for all $t\geq0$, $x\in\mathbb{R}^{d}$ a.s. That is to say
$v(t,x)(\omega)=Y(x,{\theta}_t\omega)=Y_{T-t}^{T-t,x}(\hat{\omega})$
is a stationary solution of SPDE (\ref{zz20}) w.r.t. ${\theta}$.

Therefore we have the following conclusion
\begin{thm}\label{zz41} Under Conditions {\rm(A.1)}--{\rm(A.4)}, for arbitrary $T$ and $t\in[0,T]$, let $v(t,x)\triangleq Y_{T-t}^{T-t,x}$, where $(Y_{s}^{t,x},Z_{s}^{t,x})$ is the solution of BDSDE (\ref{zz36}) with $\hat{B}_s={B}_{T-s}-{B}_T$ for all $s\geq0$. Then $v(t,x)$ is
a ``perfect" stationary stochastic viscosity solution of SPDE
(\ref{zz20}).
\end{thm}
{\bf Acknowledgements}. I would like to thank Prof. H. Z. Zhao, with
whom the main ideas were formed when I was a student in
Loughborough-Shandong Universities joint Ph.D. programme in
stochastic analysis. It is also my great pleasure to thank R.
Hudson, K. Lu, J. Ma, S. Peng and S. Tang for useful conversations. 
Meanwhile, I would like to acknowledge the financial support of the
National Basic Research Program of China (973 Program) with Grant
No. 2007CB814904.

\end{document}